\documentclass[reqno,11pt]{amsart}
\usepackage{amsmath}
\usepackage{amsthm}
\usepackage{graphicx}
\usepackage{amsfonts}
\usepackage{amssymb}
\usepackage{bm}
\usepackage{color}
\usepackage{hyperref}
\usepackage{enumerate}
\usepackage[margin=1.2in]{geometry}
\numberwithin{equation}{section}

\newtheorem{theorem}{Theorem}[section]
\newtheorem{lemma}[theorem]{Lemma}

\newtheorem{proposition}[theorem]{Proposition}
\newtheorem{corollary}[theorem]{Corollary}

\theoremstyle{definition}
\newtheorem{example}[theorem]{Example}

\newtheorem{remark}[theorem]{Remark}

\def\E{{\mathbb E}}

\def\R{{\mathbb R}}

\def\I{{\mathcal I}}
\def\H{{\mathcal H}}
\def\N{{\mathbb N}}
\def\PP{{\mathbb P}}
\def\FF{{\mathbb F}}

\def\L{{\mathcal L}}

\def\Var{\mathrm{Var}}
\def\Cov{\mathrm{Cov}}
\def\Corr{\mathrm{Corr}}

\def\CPI{C_{\mathrm{P}}}
\def\CPIA{C_{\mathrm{add}}}
\def\FC{\mathsf{FTC}}
\def\BL{d_{\mathrm{BL}}}
\def\HC{\mathsf{TC}}
\def\MMD{\mathfrak{D}}
\def\tr{\mathrm{tr}}
\DeclareMathOperator*{\avg}{avg}

\title[Approximate independence of typical marginals]{Approximate independence of typical marginals of high-dimensional probability measures}
\author{Manuel Arnese and Daniel Lacker}

\begin{document}

\begin{abstract}
For many probability measures on $\R^n$, and for $k$ sufficiently smaller than $n$, most $k$-dimensional coordinate marginals are close to product measures.
We prove this in two regimes, which notably involve no exchangeability assumptions as in prior theorems of de Finetti type.
First, for log-concave measures, the quantitative bound depends on the measure only through the operator norm of its correlation matrix, and $k$ may grow polynomially with $n$.
A key ingredient of the proof is a new strengthening of the subadditivity of relative entropy for strongly log-concave measures, which as a byproduct yields sharp bounds in the strongly log-concave case.
Second, for a broad class of measures satisfying just a variance inequality much weaker than Poincar\'e's, $k$ must stay below the sharp threshold $2\log n/\log \log n$.
Here the basis of $\R^n$ is fixed a priori, which makes this phenomenon quite different from known Gaussian approximations for low-dimensional projections such as the central limit theorem for convex bodies, although there are remarkable parallels. 
Insofar as Gaussian approximations are driven by the thin shell property, our results are driven by variance inequalities for additive functions which have appeared in the literature on the Gibbs sampler.
\end{abstract}

\maketitle

\section{Introduction}

A central feature of exchangeable interacting many-particle systems of mean field type, going often by the name of (propagation of) chaos, is that low-dimensional joint distributions are close to product measures, i.e., approximately independent. This is related to a fundamental result of Diaconis and Freedman \cite{diaconis1980finite}, which states that the $k$-dimensional marginal of any exchangeable probability measure on $\R^n$ has total variation distance at most $k^2/n$ from a mixture of iid's.  This is a finitary version of (and can be used to prove) the classical theorem of de Finetti and Hewitt-Savage characterizing infinite exchangeable sequences as mixtures of iid. When the mixing measure, which is the law of the empirical distribution, is well-concentrated (close to a point mass), it follows that the low-dimensional marginals of the original measure are close to iid.

Certain structured non-exchangeable distributions have been observed to share this feature of approximately independent low-dimensional marginals, typically due to underlying mechanisms which are quite different from the mean field situation in which interactions are averaged over all coordinates. For example, in spin systems on sparse graphs and lattices, distant vertices are approximately independent due to spatial correlation decay \cite{DemboMontanari,LackerRamananWu}, and in disordered systems such as spin-glasses, it is cancellations among the disorder variables themselves which prevent strong correlations (see \cite[Section 1.4]{Talagrand-spinglasses} for the Sherrington-Kirpatrick model and  \cite{lowe2026high} for the Hopfield model).

In this paper, we develop general criteria under which a typical (i.e., chosen uniformly at random) $k$-dimensional marginal of a probability measure on $\R^n$ is quantitatively close to the product of its marginals, for $k$ sufficiently small relative to $n$.
Precise results are stated below in two different regimes. The first is log-concave measures, for which the conclusions are the strongest: a typical $k$-dimensional marginal is close to a product measure as long as $k=o(n^{1/4}/\log^{1/2} n)$, and this can be pushed up to the sharp $k=o(\sqrt{n})$ under upper and lower bounds on the Hessian of the log-density. The second regime applies to general distributions satisfying a variance inequality  which resembles the Poincar\'e inequality but is significantly weaker. In this regime, $k$ can grow  at most like $\delta\log n / \log\log n$, for any $\delta < 2$, and this threshold $\delta=2$ is sharp.

These results are remarkably parallel to known results on a related but distinct phenomenon, wherein a typical $k$-dimensional \emph{projection} is close to Gaussian. Indeed, in the log-concave regime, the central limit theorem of Klartag \cite{klartag2007central,klartag2007power} shows that the typical $k$-dimensional projection of an $n$-dimensional isotropic log-concave measure is close to standard Gaussian, with $k$ allowed to grow polynomially in $n$. Our second regime parallels a result of Meckes \cite{meckes2012projections}, which applies to general distributions satisfying a thin shell estimate and yields a Gaussian approximation for $k = O(\delta\log n / \log\log n)$, for  any $\delta < 2$. She gives a counterexample when $\delta > 2$.

For another parallel, we recall another result of Diaconis-Freedman \cite{diaconis1987dozen}, which states that the $k$-dimensional marginal of any rotation-invariant  probability measure on $\R^n$ has total variation distance $O(k/n)$  from a scale mixture of Gaussians. This is a finitary version of (and can be used to prove)  a result attributed to Schoenberg which characterizes infinite rotation-invariant sequences as mixtures of iid Gaussians. When the mixing measure, which is the law of the Euclidean norm, is well-concentrated (known as the thin shell condition), it follows that the $k$-dimensional marginal is close to iid Gaussian.

To understand why the typical low-dimensional projections of a general, non-rotation-invariant measure should be close to a Gaussian (or mixture), one can combine this Diaconis-Freedman bound with a remarkable observation of Sudakov \cite{sudakov1978typical}: the concentration of the uniform (Haar) measure on $SO(n)$ implies that most rotations of a measure on $\R^n$ are close to its $SO(n)$-averaged version (which is, of course, rotation-invariant).
In a sense, Sudakov's observation reduces the problem of studying typical low-dimensional projections to the rotation-invariant case.

This is where the parallel breaks with our coordinate-marginal setting; there is no analogous reduction to the exchangeable (permutation-invariant) case. The uniform measure on the symmetric group $S_n$ does not concentrate nearly as strongly as that of $SO(n)$, and it is simply not the case that typical permutations of a measure on $\R^n$ must be close to the permutation-averaged version.
From this perspective, it is not a priori clear that one should expect for typical coordinate marginals to be approximately independent in much generality. 
Any such result must operate through a different mechanism, rather than concentration on the symmetric group.
The right mechanism, as we will explain, seems to be variance inequalities for additive functions. Similar inequalities have appeared in the literature on the mixing rate of the Gibbs sampler and other related Markov chains, as noteworthy consequences of stronger approximate tensorization inequalities.

\subsection*{Acknowledgements}
D.L.\ acknowledges support from an Alfred P.\ Sloan Fellowship and the NSF CAREER award DMS-2045328. ChatGPT was used to assist with some proof development and check for errors. The authors wrote the entire article and take full responsibility for its content.

\section{Main results}

For a random vector $X\sim\rho$ in $\R^n$ and a set of coordinates $v \subset [n]=\{1,\ldots,n\}$, let $\rho_v$ denote the corresponding marginal law of $X_v=(X_i)_{i \in v}$. When $v$ contains only one or two indices $\{i\}$ or $\{i,j\}$, we will use the shorthand $\rho_i$ and $\rho_{ij}$, respectively.
When $\rho$ admits a density, we will denote it by the same letter: $\rho(dx)=\rho(x)dx$.

For a function $F$ on subsets of $[n]$, we write $\avg_{|v|=k} F(v)$ to denote the average over all choices of $v$ of cardinality $k$.
The covariance matrix is defined as usual by
\[
\Cov(\rho) = \Cov(X) = \E[X^{\otimes 2}] - (\E X)^{\otimes 2},
\]
and we say that $\rho$ is isotropic if $\Cov(\rho)=I$.
Let us say that $\rho$ is \emph{coordinate-wise non-degenerate} if $\Var(X_i) > 0$ for all $i$, and in this case we define the correlation matrix 
\[
\Corr(\rho)_{ij} = \Corr(X)_{ij} = \frac{\Cov(X_i,X_j)}{\sqrt{\Var(X_i)\Var(X_j)}}.
\]
This matrix is always positive semidefinite. The total variation of a signed measure $m$ is defined as the supremum of $\int f\,dm$ over all measurable functions $|f| \le 1$.  The bounded-Lipschitz distance between two probability measures $\BL(\mu,\nu)$  is the supremum of $\int f\,d\mu - \int f\,d\nu$ over all 1-Lipschitz functions $f$ satisfying $|f| \le 1$.

\subsection{Log-concave measures}

Our first main result applies to log-concave measures. Recall that a function $\rho : \R^n \to [0,\infty)$ is called \emph{log-concave} if
\[
\rho(tx+(1-t)y) \ge \rho(x)^t\rho(y)^{1-t}, \quad \forall x,y \in \R^n, \ t \in (0,1).
\] 
Here and throughout the paper, we write $x \lesssim y$ to mean that $x \le cy$ for a universal (numerical) constant $c$.

\begin{theorem} \label{th:main}
Let $\lambda \ge 1$.
Let $\rho$ be a log-concave probability density on $\R^n$.
\begin{itemize}
\item If $\rho$ is coordinate-wise non-degenerate with $\Corr(\rho) \le \lambda I$, then 
\begin{equation}
\avg_{|v|=k}\|\rho_v-\otimes_{i \in v}\rho_i\|_{\mathrm{TV}} \lesssim \lambda^{1/3} \frac{k  \log^{1/3}n}{n^{1/6}} , \qquad 2 \le k \le n. \label{ineq:main:TV}
\end{equation}
\item If instead $\rho$ satisfies $\Cov(\rho) \le \lambda I$, then
\[
\avg_{|v|=k}\BL(\rho_v,\otimes_{i \in v}\rho_i) \lesssim  \lambda^{1/3} \frac{k^{2/3}  \log^{1/3}n}{n^{1/6}} , \qquad 2 \le k \le n.
\]
\end{itemize}
\end{theorem}

The proof is given in Section \ref{se:proof:main}.
The exponents in Theorem \ref{th:main} are almost certainly suboptimal. A noteworthy feature of this theorem is that the dimension $k$ can grow like a power of $n$, as is also the case in the version of the central limit theorem for  log-concave measures recalled below as Theorem \ref{th:KlartagCLT}.  Corollary \ref{co:main-stronger} and Proposition \ref{pr:GaussianCase} below suggest that the optimal bound might be $O(k/\sqrt{n})$, ignoring dependence on $\lambda$.

Note that the left-hand side of \eqref{ineq:main:TV} is invariant under diagonal affine transformations, i.e., pushforwards of $\rho$ by maps of the form $x \mapsto x_0+Dx$ for $x_0 \in \R^n$ and a diagonal $n \times n$ matrix $D$ with positive entries. The correlation matrix $\Corr(\rho)$ shares this invariance. Note that the largest eigenvalue of $\Corr(\rho)$ is always no more than the ratio of the largest to smallest eigenvalue of $\Cov(\rho)$.

Some dependence on $\lambda$ is necessary in both parts of Theorem \ref{th:main}. For instance, suppose $\rho$ is the centered Gaussian with covariance matrix $\Sigma = I + c(J-I)$ where $J$ is the all-ones matrix, and $0 < c < 1$. Then $\Corr(\rho)=\Cov(\rho)$ has largest eigenvalue $1+c(n-1)$, and it is straightforward to see that $\|\rho_v-\otimes_{i \in v}\rho_i\|_{\mathrm{TV}}$ is a strictly positive function of $c > 0$, which does not depend on $n$, and which depends on $v$ only through its cardinality.

Theorem \ref{th:main} applies to the uniform measure $\rho$ on any convex body $K \subset \R^n$. Letting $K_v$ denote the support of the marginal $\rho_v$, note that
\[
\big(\otimes_{i \in v}\rho_i\big)(\times_{i \in v}K_i \setminus K_v) \le \|\rho_v-\otimes_{i \in v}\rho_i\|_{\mathrm{TV}}.
\]
Hence, from Theorem \ref{th:main} it follows that most low-dimensional coordinate projections of $K$ are ``close'' to rectangles, in a measure-theoretic sense. We make no statements, however, about geometric notions of closeness such as Hausdorff distance.

The key ingredient in the proof of Theorem \ref{th:main} is a new entropy inequality which improves upon the standard subadditivity inequality for strongly log-concave measures. Let $H$ denote relative entropy (a.k.a.\ Kullback-Leibler divergence), and define
\[
\HC(\rho_v) := H(\rho_v\,|\,\otimes_{i \in v} \rho_i), \quad v \subset [n].
\]
This is known in information theory as \emph{total correlation} \cite{watanabe1960information}, or \emph{mutual information} more commonly in  the case $|v|=2$.
The classical subadditivity inequality of entropy \cite[Theorem 1]{DemboCoverThomas} implies that\footnote{In fact, because the reference measure in $\HC(\rho)=H(\rho\,|\,\otimes_{i \in [n]}\rho_i)$ is the product of the marginals of $\rho$, the constant $k/n$ in this subadditivity inequality can be sharpened to $(k-1)/(n-1)$ by a short argument of conditioning on one variable at a time, but this detail will not be important to us.}
\[
\avg_{|v|=k} \HC(\rho_v) \le \frac{k}{n} \HC(\rho).
\]

This holds generally for any probability measure $\rho$ on any measurable space.
For $\rho$ satisfying $\HC(\rho)=o(n)$, one can conclude that the left-hand side is $o(1)$ for fixed $k$, meaning that a typical low-dimensional marginal of $\rho$ is close to being equal to the product of its marginals, i.e., close to being independent.
But this general subadditivity bound is not strong enough to lead to anything like Theorem \ref{th:main},  as $\HC(\rho)$ is quite often of order $n$ rather than $o(n)$; for example, take $\rho$ to be the $(n/2)$-fold product of a non-product measure on $\R^2$.
The following theorem strengthens (for $k \ll n$) the dependence on $k/n$ from linear to quadratic, for a class of reference measures.

\begin{theorem} \label{th:independent}
Let $b \ge a > 0$, and let $\lambda=b/a$.
Let $\rho$ be a probability measure on $\R^n$ with strictly positive twice continuously differentiable density, also denoted $\rho$, which satisfies $aI \le -\nabla^2\log\rho(x)$ and $-\partial_{ii}\log\rho(x) \le b$ for all $x \in \R^n$ and $i \in [n]$.
Then, for integers $2 \le k \le n$,
\[
\avg_{|v|=k} \HC(\rho_v) \le \lambda \frac{k(k-1)}{n(n-1)}\HC(\rho).
\]
\end{theorem}

The proof is given in Section \ref{se:proof:independent}.
Using standard information-theoretic inequalities to bound $\HC(\rho)$, we deduce a stronger version of Theorem \ref{th:main} for measures $\rho$ as in Theorem \ref{th:independent}:

\begin{corollary} \label{co:main-stronger}
Let $b \ge a > 0$, and set $\lambda=b/a$. Let $\rho$ be as in Theorem \ref{th:independent}. Then $\HC(\rho)\le (n/2)\log \lambda$. In particular, for integers $2 \le k \le n$,
\begin{equation}
\avg_{|v|=k} \HC(\rho_v) \le  \frac12 \lambda \log \lambda  \,  \frac{k(k-1)}{n-1} . \label{ineq:consequence-b}
\end{equation}
\end{corollary}

The details of this corollary are given at the end of Section \ref{se:proof:independent}, along with the proof of the following proposition which records the Gaussian case, where the same dependence on $k$ and $n$ appears:

\begin{proposition} \label{pr:GaussianCase}
Let $\lambda \ge 1$.
Let $\rho$ be a nondegenerate Gaussian distribution in $\R^n$ with $\Corr(\rho) \le \lambda I$. Then, for integers $2 \le k \le n$,
\[
\avg_{|v|=k}\|\rho_v-\otimes_{i \in v}\rho_i\|_{\mathrm{TV}}^2 \le  18(\lambda-1)\frac{k(k-1)}{n-1} .
\]
\end{proposition}

\begin{remark}[Sharpness] { \ } 
\begin{enumerate}
\item The dependence on $k$ and $n$ is sharp in both Corollary \ref{co:main-stronger} and Proposition \ref{pr:GaussianCase}. To see this, let $(X_1,\ldots,X_n)\sim\rho$ be a centered Gaussian in even dimension $n=2m$, where each $X_i$ has variance 1, and $\Cov(X_i,X_{i+1})=c \in (0,1)$ for odd $i$, with all other entries of the covariance matrix vanishing. 
Since total variation does not increase when coordinates are dropped, there is a constant $\delta>0$ depending only on $c$ such that $\|\rho_v-\otimes_{i \in v}\rho_i\|_{\mathrm{TV}}^2 \ge \delta$ whenever $v$ contains $\{i,i+1\}$ for some odd $i$. A uniformly random choice of $v \subset[n]$ of size $k$ contains no such pair with probability
\[
2^k \frac{{n/2 \choose k}}{{n \choose k}} = \prod_{j=0}^{k-1} \Big(1-\frac{j}{n-j}\Big) \le e^{-k(k-1)/2n}.
\]
Hence,
\[
\avg_{|v|=k}\|\rho_v-\otimes_{i \in v}\rho_i\|_{\mathrm{TV}}^2 \ge \delta(1-e^{-k(k-1)/2n}),
\]
which is asymptotic to  $\delta k(k-1)/2n$  as $n\to\infty$ for $k = o(\sqrt{n})$. By Pinsker's inequality, example also shows that the dependence on $k$ and $n$ in Theorem \ref{th:independent} is sharp.

\item Proposition \ref{pr:GaussianCase} fails if the left-hand side is replaced by $\avg_{|v|=k}\HC(\rho_v)$, even if $\lambda-1$ on the right-hand side is replaced by an arbitrary finite-valued function of $\lambda$.  Indeed, suppose $X=(X_1,\ldots,X_n)$ is a centered Gaussian, with $(X_1,X_2)$ being independent of $(X_3,\ldots,X_n)$, with $\Var(X_1)=\Var(X_2)=1$, with $\Cov(X_1,X_2)=c \in (0,1)$, and with $X_3,\ldots,X_n$ being iid standard Gaussians.

Then $\Corr(\rho) =\Cov(\rho) \le (1+c)I \le 2I$, but $\HC(\rho_v)=-\frac12\log(1-c^2)$ if $\{1,2\} \subset v$ and so
\[
\avg_{|v|=k}\HC(\rho_v)=-\frac12\log(1-c^2)\frac{k(k-1)}{n(n-1)}.
\]
This diverges  as $c\to 1$. A lower bound away from zero is needed on $\Corr(\rho)$ in order to avoid this situation.
\item The linear dependence on $\lambda$ is essentially the best possible in Theorem \ref{th:independent}, as we illustrate in Section \ref{se:kappa-dependence}. We do not know if the dependence on $\lambda$ is sharp  in Corollary \ref{co:main-stronger}.
\end{enumerate}
\end{remark}

The proof of Theorem \ref{th:independent} relies on the following variance inequality for additive functions: For $X=(X_1,\ldots,X_n)\sim\rho$ as in Theorem \ref{th:independent}, and for square-integrable $f_1,\ldots,f_n$,
\begin{equation}
\Var\Big(\sum_{i=1}^n f_i(X_i) \Big) \le \frac{b}{a}\sum_{i=1}^n\Var(f_i(X_i)). \label{ineq:variance-main-intro}
\end{equation}
By replacing $f_i$ by $c_if_i$ and optimizing over $(c_1,\ldots,c_n) \in \R^n$, it is straightforward to check that \eqref{ineq:variance-main-intro} is equivalent to
\begin{equation}
\big\|\Cov\big((f_i(X_i))_{i=1}^n\big)\big\|_{\mathrm{op}} \le b/a, \quad \forall f_i \text{ s.t. } \E f_i(X_i)=0, \ \E f_i^2(X_i) \le 1. \label{ineq:variance-main-restated}
\end{equation}
That is, $b/a$ bounds the largest eigenvalue of the covariance matrix of any random vector of the form $(f_1(X_1),\ldots,f_n(X_n))$. Up to a normalization, this coincides with what was called the \emph{summation-based correlation coefficient} in \cite{QinWang}.

The inequality \eqref{ineq:variance-main-intro} is deduced in Lemma \ref{le:keyinequality} as a quick consequence of an approximate tensorization inequality recently discovered by Ascolani-Lavenant-Zanella in their study of the mixing time of the Gibbs sampler for log-concave measures \cite{AscolaniLavenantZanella}.
Note that \eqref{ineq:variance-main-intro} holds trivially if the constant $b/a$ is replaced by $n$. On the opposite extreme, if the $X_i$ are independent, then $b/a$ can be replaced by $1$. For $1 < b/a < n$, the inequality \eqref{ineq:variance-main-intro} encodes something about the degree of independence of $X_1,\ldots,X_n$.
This is closely related to the notion of \emph{spectral independence} codified in \cite{AnariLiuOveisGharan} for measures on the hypercube.
In our proof of Theorem \ref{th:independent}, we make key use of the fact that the inequality \eqref{ineq:variance-main-intro} holds not only for $\rho$ but also all of its conditional measures, which inherit the same strong log-concavity and smoothness parameters $a$ and $b$.

This leads first to an inequality like Theorem \ref{th:independent} but with relative Fisher information replacing relative entropy, which we integrate using the de Bruijn identity to get the announced entropy inequality. See Section \ref{se:proof:independent} for the details.

Theorem \ref{th:main} is deduced from Corollary \ref{co:main-stronger} by approximating a generic log-concave measure by one with uniform upper and lower Hessian bounds. We impose an upper Hessian bound by convolution with a Gaussian, and a lower bound using  Eldan's stochastic localization. This program is carried out in Section \ref{se:proof:main}.

It is worth mentioning that our results easily adapt to the case where each coordinate of the random vector $X=(X_1,\ldots,X_n)\sim\rho$ takes values  not in $\R$ but in some other finite-dimensional vector space. Indeed, the result of \cite{AscolaniLavenantZanella} that we rely on is stated at this level of generality. The proofs would not change much, just requiring tracking the dependence on the dimensions of these spaces.

\subsection{Beyond log-concavity}

We next move beyond the log-concave class. A different variance inequality for additive functions will serve as our key assumption. Define $\CPIA(\rho)$ to be the smallest  $c\ge 0$ for which the inequality
\begin{equation}
\Var\Big(\sum_{i=1}^nf_i(X_i)\Big) \le c\sum_{i=1}^n \|f_i'\|_{L^\infty(\R)}^2 \label{ineq:Poincare-additive}
\end{equation}
holds for all continuously differentiable functions $f_i \in L^2(\rho_i)$, $i=1,\ldots,n$. Here $X=(X_1,\ldots,X_n) \sim \rho$.
We set $\CPIA(\rho)=\infty$ if no such constant exists.
The inequality \eqref{ineq:Poincare-additive} is similar in structure to \eqref{ineq:variance-main-intro}, just involving a different squared seminorm on the right-hand side. Analogous to \eqref{ineq:variance-main-restated}, by a simple scaling argument, $\CPIA(\rho)$ can be equivalently defined as 

\begin{equation*}
\CPIA(\rho) = \sup\Big\{\big\|\Cov\big((f_i(X_i))_{i=1}^n\big)\big\|_{\mathrm{op}} : f_1,\ldots,f_n \text{ 1-Lipschitz}\Big\}.
\end{equation*}

\begin{theorem} \label{th:BLdistance}
Every probability measure $\rho$ on $\R^n$ satisfies, for integers $2 \le k \le n$,
\[
\avg_{|v|=k}\BL(\rho_v,\otimes_{i \in v}\rho_i) \lesssim  \big(1+\CPIA(\rho)^{1/3}\big)  \frac{k^{1/3}}{ (n-1)^{\frac{2}{3(k+4)}}} .
\]
The right-most fraction vanishes if $n\to\infty$ and $k \le \delta\log n/\log \log n$  for some $\delta < 2$. 
\end{theorem}

A simple argument gives a rough sense of why the inequality \eqref{ineq:Poincare-additive} should control the approximate independence of marginals, in the $k=2$ case. Take 1-Lipschitz test functions $f_1,\ldots,f_n,g_1,\ldots,g_n$ which are bounded in absolute value by 1. Applying Cauchy-Schwarz and \eqref{ineq:Poincare-additive} yields
\[
\frac{1}{n^2}\sum_{i,j=1}^n \Cov(f_i(X_i),g_j(X_j)) = \Cov\Big(\frac{1}{n}\sum_{i=1}^n f_i(X_i),\frac{1}{n}\sum_{j=1}^n g_j(X_j)\Big) \le \frac{1}{n}\CPIA(\rho).
\]
Each diagonal term $i=j$ in the double sum is $O(1)$, and we deduce that the average value of $\Cov(f_i(X_i),g_j(X_j))$ over all distinct pairs $(i,j)$ is of order $1/n$.
We may take a supremum over the choice of $f_1,\ldots,f_n,g_1,\ldots,g_n$ outside of this average over $(i,j)$, which results in a sort of approximate independence result for a typical two-dimensional marginal. But this is an extremely weak quantification of approximate independence of typical 2-dimensional marginals, unless we could somehow bring the supremum inside the average. And this argument does not extend easily from $k=2$ to $k>2$. The proof of Theorem \ref{th:BLdistance}, given in Section \ref{se:proof:BLdistance}, follows a completely different route.

As we will explain in Section \ref{se:intro-coordinatemarginals}, the analogy with Gaussian approximations of typical projections might lead one to guess that $\CPIA(\rho)$ could be replaced by a parameter measuring the concentration of the empirical measure. Example \ref{ex:empiricalmeasureconcentrates} shows that this is not possible.

Log-concavity implies $\CPIA(\rho) <\infty$, because $\CPIA(\rho)$ is no greater than the usual Poincar\'e constant. 
The conclusion of Theorem \ref{th:BLdistance} is significantly weaker than that of Theorem \ref{th:main} in the dependence on $n$ and $k$; indeed, the latter is informative even when $k$ grows polynomially with $n$. In fact, the threshold $k\sim 2\log n/ \log\log n$ in Theorem \ref{th:BLdistance} is sharp:

\begin{proposition} \label{pr:sharpness-log}
For every $\delta > 2$ there exists $N_\delta >0$ such that for every $n \ge N_\delta $ there exists a probability measure $\rho$ on $\R^n$ such that $\CPIA(\rho) = 1$ and  such that,  for every integer $k$ with $n \ge k \ge \delta\log n/ \log\log n$, 
\[
\min_{|v|=k}\BL(\rho_v,\otimes_{i \in v}\rho_i) \ge 1/2.
\]
\end{proposition}

The proof is given in Section \ref{se:proof:pr:sharpness-log}.

Let us  describe some classes of examples of  measures satisfying \eqref{ineq:Poincare-additive} and its relation to  other inequalities.
Its closest relative is perhaps the $(2,\infty)$-Poincar\'e inequality restricted to additively separable functions. Recall that $\rho$ satisfies  a $(p,q)$-Poincar\'e inequality if
\begin{equation}
\Big\| f - \int f\,d\rho \Big\|_{L^p(\rho)} \le C_{p,q} \| |\nabla f| \|_{L^q(\rho)}, \label{ineq:Poincare-pq}
\end{equation}
for all continuously differentiable functions $f \in L^p(\rho)$.
Note that the usual Poincar\'e inequality corresponds to $p=q=2$, and we write $\CPI(\rho)$ for the square of the optimal constant $C_{2,2}$.
In the case $(p,q)=(2,\infty)$, applying \eqref{ineq:Poincare-pq} to a function of the form $f(x)=\sum_{i=1}^nf_i(x_i)$ shows that
\[
\|\Cov(\rho)\|_{\mathrm{op}} \le \CPIA(\rho) \le C_{2,\infty}^2 \le C_{2,2}^2=\CPI(\rho).
\]
The left-most inequality here follows by testing \eqref{ineq:Poincare-additive} against (bounded approximations of) linear functions $f_i$.
For log-concave measures $\rho$, the recently resolved Kannan-Lova\'sz-Simonovits conjecture \cite{BizeulKlartagLehec-KLS} asserts  that $\CPI(\rho) \lesssim \|\Cov(\rho)\|_{\mathrm{op}}$, and thus the above chain of inequalities collapses to equivalence up to universal constants.
Moreover, a seminal result of E.\ Milman \cite[Theorem 2.4]{milman2009role} shows that the constants $C_{p,q}$ are all quantitatively equivalent for log-concave measures, in the sense that $C_{p',q'} \lesssim p'C_{p,q}$. 

In the non-log-concave case, the  inequality \eqref{ineq:Poincare-additive} can often be checked much more easily than any $(p,q)$-Poincar\'e inequality. Indeed, we have
\[
\Var(f(X)) \le \E[(f(X)-f(\E X))^2] \le 
\| |\nabla f| \|_{L^\infty(\R^n)}^2\sum_{i=1}^n\Var(X_i).
\]
This yields $\CPIA(\rho) \le \sum_{i=1}^n\Var(X_i)$, which is at most  $D^2/2$ if  $\rho$ has compact support with diameter $D$.
Notably, this is valid even when the support of $\rho$ is disconnected, in which case the standard Poincar\'e inequality fails. (Note that this is possible because we use the $L^\infty(\R)$ norm in \eqref{ineq:Poincare-additive} instead of $L^\infty(\rho_i)$.)

For another class of examples, if $X\sim \rho$ is supported on the discrete hypercube $\{0,1\}^n$, then $\|\Cov(\rho)\|_{\mathrm{op}}=\CPIA(\rho)$. Indeed, the inequality $\le$ was mentioned above, and it becomes equality here because all functions on $\{0,1\}$ are restrictions of affine functions on $\R$.

Thus, Theorem \ref{th:BLdistance} is informative for any measure on the hypercube satisfying $\|\Cov(\rho)\|_{\mathrm{op}} = o(n^{2/(k+4)})$. This covers many noteworthy examples, such as Ising measures of the form
\[
\rho(x) \propto \exp\big(-\beta \langle x,Jx\rangle/2\big), \qquad x \in \{-1,1\}^n,
\]
where $J$ is a symmetric matrix. If the difference $\delta$ between the largest and smallest eigenvalues of $J$ is strictly smaller than $1/\beta$, then $\|\Cov(\rho)\|_{\mathrm{op}}$ is bounded by a constant depending only on $\beta\delta$; this follows from the log-Sobolev inequality of \cite[Theorem 1]{bauerschmidt2019very}.
For the Sherrington-Kirpatrick (SK) model, where $J$ is a GOE matrix with upper triangle consisting of independent centered Gaussians of variance $1/n$, and the largest and smallest eigenvalues converge to $\pm 2$, the preceding statement applies if $\beta < 1/4$. We can improve this using sharper results specific to the SK model: If $\beta<1$, it is known that $\|\Cov(\rho)\|_{\mathrm{op}} \to 1/(1-\beta)^2$ in probability as $n\to\infty$ \cite{brennecke2024two}, and $\E\|\Cov(\rho)\|_{\mathrm{op}}$ is bounded by a constant depending only on $\beta$ \cite{elAlaoui2024bounds}; boundedness of $\|\Cov(\rho)\|_{\mathrm{op}}$ in probability was extended to the case with external field in \cite{brennecke2}. In all of these cases, we deduce from Theorem \ref{th:BLdistance} that most $k$-dimensional coordinate marginals are close to product measures, for $k$ smaller than $2\log n / \log \log n$, and this should be understood in the quenched sense for the SK model.
It should be noted that in the specific models discussed in this paragraph there may be direct methods that yield sharper bounds than our Theorem \ref{th:BLdistance}; for the SK model, see for instance \cite[Theorem 1.4.15]{Talagrand-spinglasses}.

More examples of measures with $\CPIA(\rho) < \infty$ can be produced using a tensorization principle. If $\rho$ is pairwise independent, in the sense that $\rho_{ij}=\rho_i\otimes\rho_j$ for all distinct $i$ and $j$, then $\CPIA(\rho) = \max_i\CPIA(\rho_i)$. This observation is used and spelled out in more detail in the proof of Proposition \ref{pr:sharpness-log}.

\subsection{Typical coordinate marginals versus typical projections}

Here we recall some known results about typical \emph{projections}, with rotations allowed, and we explore the similarities and differences with our study of coordinate-marginals.
For a random vector $X\sim \rho$ in $\R^n$ and a subspace $E \subset \R^n$, let $\rho_E$ be the law of the projection of $X$ onto $E$.
Let $\mathsf{Gr}(k,n)$ denote the Grassmannian manifold of $k$-dimensional subspaces of $\R^n$, equipped with its uniform probability measure (which can be described concretely as the law of the span of $k$ iid standard $n$-dimensional Gaussian vectors). We write $dE$ for integration with respect to this uniform measure. Let $\gamma_E$ denote the standard Gaussian on $E$.

The following central limit theorem for log-concave measures is due to Klartag:

\begin{theorem}[Theorem 1.1 of \cite{klartag2007power}]\!\!\!\footnote{Theorem 1.1 of \cite{klartag2007power} is stated for uniform measures on convex bodies, but it is explained therein that the result is valid also for general log-concave distributions, and the proof is written at that level of generality. Moreover, the statement in \cite{klartag2007power} is a high-probability bound for a uniform choice from $\mathsf{Gr}(k,n)$, and integrating yields the version we have stated here.} \label{th:KlartagCLT}
There exist universal constants $c,\delta>0$ such that the following holds. Let $\rho$ be an  isotropic log-concave probability measure on $\R^n$ with mean zero. 
For integers $1\le k\le cn^\delta$,

\[
\int_{\mathsf{Gr}(k,n)}\|\rho_E - \gamma_E\|_{\mathrm{TV}}\,dE \le c \,n^{-\delta}.
\]
\end{theorem}

In words, $\rho_E$ is close to standard Gaussian for most choices of $E$.
As with our Theorem \ref{th:main}, $k$ can grow polynomially in $n$.

In general, little is known about precisely which subspaces $E$ this statement applies to, except when $k=1$ in special cases such as unconditional measures \cite[Theorem 1.1]{klartag2007central}.
Obviously the bound of Theorem \ref{th:KlartagCLT} cannot apply in a pre-specified basis (that is, to an arbitrary choice of $E$), as evidenced by product measures.

It is well understood that  log-concavity is not essential to this phenomenon of typical projections being Gaussian. Rather, the fundamental ingredient is a thin shell estimate, quantifying how well concentrated $|X|$ is. The recent textbook \cite{BobkovChistyakovGotze} documents state-of-the-art quantitative bounds on various distances between typical one-dimensional projections and the standard Gaussian, under a variety of concentration assumptions of different strengths; see also the recent articles \cite{klartag2026revisiting,XHu} giving sharp bounds in different metrics.
These works deal with the $k=1$ case.
A close parallel to our Theorem \ref{th:BLdistance} is the following (isotropic case of a) theorem of E.\ Meckes: 

\begin{theorem}[Theorem 1.3 of \cite{meckes2012projections}] \label{th:MeckesCLT}
Let $\rho$ be an isotropic probability measure on $\R^n$ with mean zero.
Then, for $1\le k\le n$,
\begin{equation*}
\int d_{\mathrm{BL}} \big(\rho_E,\gamma_E\big)\,dE
\lesssim \frac{k+\log n}{k^{2/3}n^{2/(3k+4)}}+\frac{k+\sqrt{k \Var(|X|^2)}}{n}.
\end{equation*}
\end{theorem}

The recently resolved thin shell conjecture \cite{klartag2025thin}, or its sharp form \cite{ChenKlartag}, gives $\Var(|X|^2)\le 8n$ for centered isotropic log-concave $X$.
Even without log-concavity, when this thin-shell parameter $\Var(|X|^2)$ is, say, $O(n^p)$ for some $p < 2$, Theorem \ref{th:MeckesCLT} yields approximate Gaussianity for typical $k$-dimensional projections when $k \le \delta\log n/\log \log n$ for some $\delta < 2$. Moreover, Meckes constructs an example showing that in general $k$ can be no larger than $\delta\log n / \log\log n$ for $\delta > 2$. Remarkably, the critical scale $\delta=2$ is exactly the same in her work and ours. Her example is the uniform distribution over the vertices of the cross-polytope $\{\pm e_1,\ldots,\pm e_n\}$. This example does not work for our purposes; both $\rho_v$ and $\otimes_{i \in v}\rho_i$ have a point mass at zero of size at least $1-\frac{k}{n}$ for any $v \subset [n]$ with $|v|=k$, so their TV distance is at most $2k/n$. Our construction for Proposition \ref{pr:sharpness-log} starts from a discrete measure of much smaller support, built using ideas from combinatorial design theory, then smooths it into a continuous measure with uniform two-dimensional marginals.

Despite the close parallels between our results on \emph{typical coordinate-marginals} (Theorems \ref{th:main} and \ref{th:BLdistance}) and those on  \emph{typical projections} (Theorems \ref{th:KlartagCLT} and \ref{th:MeckesCLT} of Klartag and Meckes), the mechanisms behind them turn out to be quite different. Let us begin by summarizing the projection setting.

\subsubsection{Typical projections}
The first step is an observation going back to Sudakov \cite{sudakov1978typical}, which is essentially that the strong concentration of the Haar measure on the special orthogonal group $SO(n)$ implies that most rotations of $\rho$ are close to its averaged version. Precisely, a Haar-random matrix $Q$ in $SO(n)$ satisfies
\[
\PP(g(Q) - \E g(Q) > t) \le \exp\big(-(n-2)t^2/24\big), \quad t > 0,
\]
for any function $g$ on $SO(n)$ which is 1-Lipschitz with respect to the Hilbert-Schmidt norm \cite[Theorem 5.17]{Meckes-HaarBook}. For $X\sim \rho$ being centered and isotropic, let $\rho_A$ denote the law of $AX$ for each (non-random) $A \in SO(n)$, and let $\overline\rho=\E \rho_Q$ denote the averaged (rotation-invariant) measure.
For $f : \R^n \to \R$ being 1-Lipschitz, the function $A \mapsto \int f\,d\rho_A$ is readily checked to be 1-Lipschitz in the Hilbert-Schmidt norm, and we deduce that there is a set $S_{t,f} \subset SO(n)$ of Haar measure at least $1-2\exp\big(-(n-2)t^2/24\big)$ such that, for each $t > 0$ and $A \in S_{t,f}$, we have $|\int f\,d\rho_A - \int f\,d\overline\rho| \le t$. As $n$ grows, choosing such an $A$ becomes overwhelmingly likely. The proof of Meckes's  Theorem \ref{th:MeckesCLT} above combines this approach with a careful chaining argument. 

This measure concentration argument allows us to reduce the study of typical projections to the study of the rotation-invariant measure $\overline\rho$. Diaconis and Freedman showed in \cite{diaconis1987dozen} that the low-dimensional marginals of a rotation-invariant measure on $\R^n$ must be close to mixtures of Gaussians. More precisely, the $k$-dimensional marginal of $\overline{\rho}$ has total variation distance at most $O(k/n)$ from $\int_{\R_+} N_k(0,rI) \eta(dr)$, for $\eta$ being the law of $|X|^2/n$ and $N_k(0,rI)$ being the centered Gaussian in $\R^k$ with covariance matrix $rI$. This is shown by quantifying the classical fact that the $k$-dimensional marginals of the uniform measure on $\sqrt{n}\mathbb{S}^{n-1}$ are close to $N_k(0,I)$, for $k\ll n$.
Finally, if the measure $\eta$ (which has mean 1 in the isotropic case) is close to the Dirac $\delta_1$, the aforementioned Gaussian mixture must in fact be close to $N_k(0,I)$. This concentration of $\eta$ 
can be quantified using the thin shell bounds discussed above. 

\subsubsection{Typical coordinate marginals} \label{se:intro-coordinatemarginals}
The previous paragraph has parallels in our setting of coordinate marginals. Diaconis and Freedman showed in a different paper \cite[Theorem 13]{diaconis1980finite} (see also Stam \cite{stam1978distance})  that any exchangeable (permutation-invariant) probability measure $\rho$ on an $n$-fold product space is close, in terms of low-dimensional marginals, to a mixture of product measures. Precisely,
\begin{equation}
\bigg\|\rho_{\{1,\ldots,k\}} - \int_{\R^n} \Big(\frac{1}{n}\sum_{i=1}^n\delta_{x_i}\Big)^{\otimes k}\,\rho(dx)\bigg\|_{\mathrm{TV}} \le \frac{k(k-1)}{n}. \label{ineq:DiaconisFreedman}
\end{equation}
That is, ``mixture of Gaussians'' in the rotation-invariant setting is replaced by ``mixture of iid'' in the permutation-invariant setting. The thin shell property in the rotation setting was the concentration of the mixing measure, which was the law of $|X|^2/n$. In the permutation setting, the mixing measure is the law of the empirical measure $L_n=(1/n)\sum_{i=1}^n \delta_{X_i}$, which is now an infinite-dimensional object; its mean is the one-dimensional marginal $\rho_1$. Nonetheless, if $L_n$ concentrates in a suitable sense, we should expect  $\rho_{\{1,\ldots,k\}}$ to be close to $\rho_1^{\otimes k}$.  
Note that \eqref{ineq:Poincare-additive} implies but is much stronger than a form of concentration of $L_n$: For a 1-Lipschitz function $g$,
\[
\Var\Big( \int g\,dL_n\Big) \le \CPIA(\rho)/n.
\]
Exchangeability is critical to this discussion so far.

To complete the analogy, we would need to remove the assumption of exchangeability. This would require an analog of Sudakov's aforementioned observation, stating that most permutations of the measure are close to the permutation-averaged version.
But this fails dramatically, essentially because the rate of measure concentration on the permutation group $S_n$ is much slower than on $SO(n)$.
Here is a simple class of examples of $\rho$ whose coordinate marginals are not close to product measures, despite the empirical measure being well-concentrated. In particular, a direct analog of a \emph{thin shell} assumption in our context would be some notion of concentration of the empirical measure, but no such notion would suffice for Theorem \ref{th:BLdistance}.

\begin{example} \label{ex:empiricalmeasureconcentrates}
Let $n=2m$ be even, and take two distinct probability measures $\mu$ and $\nu$ on $\R$. Consider the mixture
\[
\rho = \frac12(\mu^{\otimes m} \otimes \nu^{\otimes m} + \nu^{\otimes m} \otimes \mu^{\otimes m}).
\]
In other words, take a random vector of independent coordinates such that the first $m$ are distributed according to $\mu$ and the last $m$ according to $\nu$; 
then, with probability 1/2, reverse the order of this random vector. The 1-dimensional marginals are all the same, $\rho_i=(\mu+\nu)/2$, whereas the pair marginals depend on whether or not $i$ and $j$ lie in the same block of $m$ coordinates; $\rho_{ij} = (\mu^{\otimes 2} + \nu^{\otimes 2})/2$ if both $i,j \le m$ or both $i,j > m$, and 
$\rho_{ij} = (\mu\otimes\nu+\nu\otimes\mu)/2$ otherwise.
The distance  between $\rho_{ij}$ and $\rho_i \otimes \rho_j$ is thus bounded away from zero uniformly in $i$, $j$, and $n$.
Note that the empirical measure of this random vector converges to $(\mu+\nu)/2$ by the law of large numbers.

Theorem \ref{th:BLdistance} is uninformative for this example because the parameter $\CPIA(\rho)$ here is at least of order $n$ if $\mu$ and $\nu$ are distinct. To see this, fix a bounded 1-Lipschitz function $g$ on $\R$, and set $f_i=g$ for $i \le m$ and $f_i=-g$ for $i > m$. Let $Y_i$ be iid $\sim \mu$ and $Z_i \sim$ iid $\nu$, and let $U$ be an independent coin flip. If $U$ is heads, put $X=(Y_1,\ldots,Y_m,Z_1,\ldots,Z_m)$, and otherwise put $X=(Z_1,\ldots,Z_m,Y_1,\ldots,Y_m)$. Then $X \sim \rho$, and
\begin{align*}
\Var\Big(\sum_{i=1}^n f_i(X_i) \Big) \ge \Var\Big(\E\Big[\sum_{i=1}^n f_i(X_i) \,\Big|\,U\Big]\Big) =  m^2\Big(\int g\,d(\mu-\nu)\Big)^2 =: m^2\delta_g^2.
\end{align*}

For any bounded 1-Lipschitz $g$ such that $\delta_g>0$, the left-hand side is at most $n\CPIA(\rho)$. Using $m=n/2$, it follows that $\CPIA(\rho) \ge n\delta_g^2/4$.
\end{example}

\begin{remark} \label{re:finiteexchangeablerates}
To comment briefly on the $O(k^2/n)$ bound for exchangeable measures in \eqref{ineq:DiaconisFreedman}:  It is sharp in general \cite[Section 4]{diaconis1980finite}, achieved by comparing sampling with and without replacement (the birthday paradox), but it can be improved \cite[Theorem 3]{diaconis1980finite} to $|S|k/n$ if $\rho$ is supported on $S^n$ for a finite subset $S$ of $\R$. Aside from finite state spaces, there are two other general settings we know of in which the approximation by (mixtures of) product measures holds at the improved order of $k/n$ instead of $k^2/n$. 
One is the class of noisy permutation mixtures described in \cite[Theorem 1.4]{han2026approximate}. The other is a class of mean field interacting particle systems \cite{lacker2023hierarchies,lacker2022quantitative}, which will be discussed more in the next section. Both deal with exchangeable measures. Our Corollary \ref{co:LSI+Poincare} below gives rise to a broad class of examples which need not be exchangeable.
\end{remark}

\subsection{Relation to mean field approximation}

The notion of approximate independence studied here is related to certain notions of \emph{mean field approximations}. Note first that
\begin{equation}
\HC(\rho) = \inf_{(\mu_i)_{i \in [n]}}  H(\rho\,|\,\otimes_{i \in [n]}\mu_i ), \label{eq:HCvariational}
\end{equation}
the infimum being over product measures $\otimes_{i \in [n]}\mu_i$; the optimizer is exactly $\otimes_{i \in [n]}\rho_i$. This follows from the Pythagorean identity
\begin{equation}
H(\rho\,|\,\otimes_{i \in [n]}\mu_i ) = \HC(\rho) + \sum_{i \in [n]}H(\rho_i\,|\,\mu_i). \label{eq:pythagoras}
\end{equation}
This is worth contrasting with the defining minimization problem in meanfield variational inference \cite{blei2017variational} in Bayesian statistics, which is often applied when $\rho$ itself (or its marginals) are difficult to compute with directly. Instead, one approximates $\rho$ by the minimizer of
\[
\inf_{(\mu_i)_{i \in [n]}} H(\otimes_{i \in [n]}\mu_i  \,|\,\rho),
\]
where the arguments are reversed compared to \eqref{eq:HCvariational}.

In the study of mean field dynamics, Kac's (molecular) chaos is a common asymptotic notion of approximate independence \cite{sznitman2006topics}. Suppose for each $n$ that $X^n=(X^n_1,\ldots,X^n_n)$ is an exchangeable random vector, and suppose $\mu$ is a fixed probability measure on $\R$. We say that the triangular array $(X^n_i)$ is $\mu$-chaotic if the empirical measure $\frac{1}{n}\sum_{i=1}^n \delta_{X^n_i}$ converges weakly in probability to $\mu$ as  $n\to\infty$. This is equivalent to the statement that the joint law of $(X^n_1,\ldots,X^n_k)$ converges weakly to $\mu^{\otimes k}$ as $n\to\infty$ for each $k \in \N$; see \cite[Proposition 2.2]{sznitman2006topics}.
The definition of $\mu$-chaotic combines two notions that we prefer to keep separate. These two notions are most easily disentangled by means of the entropy decomposition \eqref{eq:pythagoras}. Letting $\rho^n_k$ denote the law of $(X^n_1,\ldots,X^n_k)$,
\begin{equation}
H(\rho^n_k\,|\,\mu^{\otimes k}) = \HC(\rho^n_k) + k H(\rho^n_1\,|\,\mu). \label{intro:exchangeable-entropy-decomp}
\end{equation}
The first term quantifies approximate independence, while the second quantifies how far the one-dimensional marginals are from the target measure $\mu$.
These two quantities are quite different in nature.
To be $\mu$-chaotic, at least if measured by relative entropy, requires both quantities to be small. In this paper we focus only on the first quantity, which is intrinsic to $\rho^n$. In contrast, controlling the second quantity $H(\rho^n_1\,|\,\mu)$ requires extrinsic information, about the relationship between $\rho^n$ and $\mu$. In the literature on mean field limits, this additional information is often that the marginal $\rho^n_1$  approximately satisfies some fixed-point equation (or PDE) which $\mu$ satisfies exactly  and which enjoys some stability properties.

\subsection{A quantitative mechanism for sharp mean field approximations}

In our results so far we have not claimed optimal dependence on $k$ or $n$, except for the $k^2/n$ bound of Corollary \ref{co:main-stronger} which cannot be improved in general. However, in Remark \ref{re:finiteexchangeablerates} we discussed prior work that improves on this $k^2/n$ bound for certain  exchangeable measures. In this section we record a simple alternative to Corollary \ref{co:main-stronger}, which  retains information about the measure $\rho$, and we explain how it leads to improvements over $k^2/n$ for certain classes of measures. 
It seems to capture, in a unified way, a wide range of different decorrelation mechanisms. We explain three such mechanisms below for which Corollary \ref{co:LSI+Poincare} gives sharp quantitative results: dense, sparse, and disordered interactions.

\begin{corollary} \label{co:LSI+Poincare}
Under the assumptions of Theorem \ref{th:independent},
\[
\avg_{|v|=k}\HC(\rho_v) \le \frac{b}{a^3}\frac{k(k-1)}{n(n-1)} \sum_{1 \le i < j \le n}\int_{\R^n} |\partial_{ij}\log \rho|^2\,d\rho.
\]
\end{corollary}

The short proof is in  Section \ref{se:proof:independent}.
Note that the diagonal  $i=j$ is excluded from the summation.
The $k=n$ case of Corollary \ref{co:LSI+Poincare} resembles an estimate from  \cite{LackerMukherjeeYeung} on the mean field variational inference error, which showed that if $-\nabla^2\log\rho \ge aI$ then
\[
\inf_{(\mu_i)_{i \in [n]}} H(\otimes_{i \in [n]}\mu_i  \,|\,\rho) \le \frac{1}{a^2}\sum_{1 \le i < j \le n} \int |\partial_{ij}\log\rho|^2\,d\mu^*,
\]
where $\mu^*=\otimes_{i \in [n]}\mu^*_i$ is the (unique) optimizer of the left-hand side. It was noted in \cite[Inequality (1.16)]{LackerMukherjeeYeung} that, when combined with the Talagrand inequality, this implies
\begin{equation}
\avg_{|v|=k} W_2^2(\mu^*_v, \rho_v) \le \frac{4k}{a^3n}\sum_{1 \le i < j \le n} \int |\partial_{ij}\log\rho|^2\,d\mu^*, \label{ineq:LMY-W2}
\end{equation}
where $W_2$ is the quadratic Wasserstein distance.
If the double sum is of order $o(n)$, this yields a quantitative estimate of approximate independence for typical low-dimensional marginals of $\rho$, similar in spirit to our main results. However, it is quite often the case that the double sum is of order $n$, in which case \eqref{ineq:LMY-W2} is uninformative, while Corollary \ref{co:LSI+Poincare} is still informative.

A common class of examples are the Gibbs measures with pairwise interactions:
\begin{equation}
\rho(dx) \propto \exp\bigg( -\sum_{i=1}^n V(x_i) - \sum_{1 \le i \leq j \le n} J_{ij}  W(x_i-x_j)\bigg). \label{eq:Gibbs}
\end{equation}
Assume that $W$ is even for simplicity, and that $J$ is a symmetric matrix with nonnegative entries such that $\max_i\sum_{j \neq i} J_{ij} \le 1$. If $0 < a \le V'' \le b_V$ and  $0 \le W'' \le b_W$, then it is straightforward to check that $aI \le -\nabla^2\log \rho \le bI$ with $b=b_V + 2b_W$. Computing $\partial_{ij}\log\rho(x) = J_{ij}W''(x_i-x_j)$ for $i \neq j$,  Corollary \ref{co:LSI+Poincare} implies
\[
\avg_{|v|=k}\HC(\rho_v) \le (b/a)^3 (k/n)^2 \tr(J^2).
\]
Natural examples, discussed also in \cite[Section 2.1]{LackerMukherjeeYeung}, arise when $J=A/d$ for $A$ being the adjacency matrix of a $d$-regular graph on vertex set $[n]$. Then $\tr(J^2)=n/d$, and the bound becomes $(b/a)^3 k^2/dn$. This is informative in the very different \emph{dense} and \emph{sparse} regimes, where $d$ is divergent or bounded, respectively.

The densest regime, where $d=cn$ for some positive constant $c$, leads to a class of  mean field Gibbs measures that has been studied in many prior papers. The bound $(k/n)^2$ matches the sharp rate of convergence to the mean field limit in this convex regime; this was shown in \cite{RenWang2025}, refining the prior work of \cite{lacker2022quantitative} which initially identified this sharp rate but only for $b_W < a$. Let us stress again that our results only address the  $\HC$ term in \eqref{intro:exchangeable-entropy-decomp} and thus  do not directly pertain to the mean field limit. In this case decorrelation happens by a law of large numbers, where many weak interactions are aggregated  as $n$ grows large. Regimes of intermediate denseness $1 \ll d \ll n$ behave similarly and were studied for interacting diffusion dynamics in \cite{LackerYeungZhou}; the uniform-in-time estimates therein apply immediately to the invariant measures, which take the same form as \eqref{eq:Gibbs}.

At the opposite extreme is the sparse regime, where $d$ stays bounded as $n\to\infty$. 
The approximate independence of low-dimensional marginals arises here for entirely different reasons. The approximate independence comes from the fact that correlations decay  with respect to graph distance, and two uniformly chosen vertices tend to be far from each other. The rate becomes $k^2/n$, which is the worst case found in Corollary \ref{co:main-stronger}. This matches the birthday paradox heuristics: For instance, if $k$ vertices are selected uniformly without replacement from the $n$-vertex cycle graph, then  the minimal graph distance between pairs of vertices is of order $\max(1,n/k^2)$, which diverges precisely when $k = o(\sqrt{n})$.
Note that the approximate independence bound of Corollary \ref{co:main-stronger} does not require that we understand or formalize this notion of correlation decay.

Continuous spin systems provide another interesting regime where Corollary \ref{co:LSI+Poincare} yields bounds which are likely sharp. Let $J=AA^\top$ where $A$ is an $n\times n$ random matrix with iid sub-gaussian entries having mean $0$ and variance $1/n$. Consider the following random measure, which appears in applications in Bayesian statistics \cite{mukherjee2026meanfieldapproachempirical,fan2025dynamical,fan2025dynamicalmeanfieldanalysisadaptive}:
\[\rho^J(dx)\propto \exp\Big(-\sum_{i=1}^n V(x_i)-\sum_{1\leq i\leq j\leq n}J_{ij}x_ix_j\Big)\,dx.\]
Assuming that $0<a \le V''\leq b$, Corollary \ref{co:LSI+Poincare} yields the quenched estimate
\begin{align*}
    \avg_{|v|=k}\HC(\rho^J_v)\leq \frac{k^2}{n^2} \frac{b + 2\max_i J_{ii}}{a^3} \sum_{1\leq i<j\leq n}J_{ij}^2\leq \frac{k^2}{n} \frac{b+2\|A\|^2_{\text{op}}}{2a^3}  \|A\|_{\textup{op}}^4. 
\end{align*}
Classical random matrix results imply that $\|A\|_\text{op}=O(1)$ with probability at least $1-2e^{-cn}$ for some $c>0$, yielding the high probability estimate $\avg_{|v|=k}\HC(\rho^J_v)\lesssim k^2/n$. Taking expectations, we obtain for every $v$ of size $k$ that $\E \HC(\rho_v^J)\lesssim k^2/n$, with the expectation applying to the randomness of $J$. In this regime, decorrelation is caused by the self-averaging of the disorder matrix $J$.

\section{Proof of Theorem \ref{th:independent} and its corollaries} \label{se:proof:independent}

A key input to our analysis in the log-concave setting is a variance inequality proven recently by Ascolani-Lavenant-Zanella \cite{AscolaniLavenantZanella} in their study of the Gibbs sampler for strongly log-concave distributions.
A main result of theirs is an approximate tensorization inequality for entropy, Theorem 3.1 in their paper: For  any probability measure $\mu$ on $\R^n$,
\[
\avg_{|v|=n-1} H(\mu_v\,|\,\rho_v) \le \Big(1-\frac{a}{bn} \Big)H(\mu\,|\,\rho).
\]
Corollary 3.7 of \cite{AscolaniLavenantZanella} notes that the linearization of this inequality is 
\begin{equation}
\Var(g(X)) \le \frac{b}{a}\sum_{i=1}^n \E \Var(g(X)\,|\,  (X_j)_{j \neq i}) . \label{ineq:ALZ-key}
\end{equation}
In fact, the assumptions of \cite{AscolaniLavenantZanella} are more general than those of our Theorem \ref{th:independent}.
The main inequality we will need is \eqref{ineq:Manuel} below, with $\lambda=b/a$, but we give a bit more information than necessary in the following lemma in an effort to clarify the nature of the key inequality.

\begin{lemma} \label{le:keyinequality}
Let $\lambda \ge 0$. Let $X\sim \rho$ be a random vector in $\R^n$.
Then 
\begin{equation}
\sum_{i=1}^n \Var\big(\E[g(X)\,|\,X_i]\big) \le \lambda\Var(g(X)), \qquad \forall g \in L^2(\rho) \label{ineq:Manuel}
\end{equation}
 if and only if  
\begin{equation}
\Var\Big(\sum_{i=1}^n f_i(X_i) \Big) \le \lambda\sum_{i=1}^n\Var(f_i(X_i)), \quad \forall f_i \in L^2(\rho_i), \ i=1,\ldots,n. \label{ineq:variance-main}
\end{equation}
Moreover, if $X\sim \rho$ for $\rho$ satisfying the assumptions of Theorem \ref{th:independent}, then both inequalities hold with $\lambda=b/a$.
\end{lemma}
\begin{proof}
Define the Hilbert space $\L:= \oplus_{i=1}^n L^2_0(\rho_i)$, where the subscript $0$ under $L^2$ indicates mean zero.
Define an operator $S : \L \to L^2_0(\rho)$ by $S(\oplus_i f_i)(x):= \sum_i f_i(x_i)$ for $x \in \R^n$.
Defining also the projections $P_i : L^2(\rho) \to L^2(\rho_i)$ by $P_ig(X_i)=\E[g(X)|X_i]$, the adjoint of $S$ can be written as 
\[
S^* : L_0^2(\rho) \to \L, \qquad S^* g=\oplus_i P_ig .
\]
Indeed, for $\oplus_i f_i \in \L$ and $g \in L^2_0(\rho)$,
\begin{align*}
\langle S(\oplus_i f_i),g\rangle_{L^2(\rho)} &= \E\bigg[ g(X)\sum_{i=1}^n f_i(X_i)\bigg] = \E\sum_{i=1}^n f_i(X_i)P_ig(X_i) = \langle \oplus_i f_i, \oplus_i P_ig\rangle_{\L}.
\end{align*}
The inequality \eqref{ineq:variance-main} can be written as
\begin{align*}
\big\|S(\oplus_i f_i)\big\|_{L^2(\rho)}^2 \le \lambda\big\|\oplus_i f_i \big\|_{\L}^2, \qquad \forall \oplus_i f_i \in \L,
\end{align*}
which is precisely a bound on the operator norm, $\|S\|_{\L \to L^2_0(\rho)}^2 \le \lambda$. The adjoint has the same norm, so the preceding inequality is equivalent to $\|S^*\|_{L_0^2(\rho)\to \L}^2 \le \lambda$, which can be written as 
\[
\big\|\oplus_i P_ig \big\|_{\L}^2 \le \lambda \|g\|_{L^2_0(\rho)}^2,
\]
which is exactly \eqref{ineq:Manuel}. This proves the claimed equivalence.

Finally, assume now that $X\sim \rho$  for $\rho$ as in Theorem \ref{th:independent}.
Let $g(X)=\sum_{i=1}^n f_i(X_i)$, and note that
\[
\Var(g(X)\,|\,(X_j)_{j \neq i})=\Var(f_i(X_i)\,|\,(X_j)_{j \neq i}).
\]
Apply the inequality \eqref{ineq:ALZ-key} followed by the law of total variance to get
\begin{align*}
\Var\Big(\sum_{i=1}^n f_i(X_i)\Big) \le \frac{b}{a}\sum_{i=1}^n  \E\Var(f_i(X_i)\,|\,(X_j)_{j \neq i}) \le \frac{b}{a}\sum_{i=1}^n  \E\Var(f_i(X_i)), 
\end{align*}
which proves that \eqref{ineq:variance-main} holds with $\lambda=b/a$.
\end{proof}

\begin{proposition} \label{pr:keyvariancebound}
Let $X\sim \rho$ with $\rho$ as in Theorem \ref{th:independent}.
For any $g_1,\ldots,g_n \in L^2(\rho)$ and $2 \le k \le n$,
\begin{align*}
\avg_{|v|=k}\sum_{i \in v}\E\big|\E[g_i(X)|X_v]-\E[g_i(X)|X_i]\big|^2 \le \frac{b}{a}\frac{k(k-1)}{n(n-1)}\sum_{i=1}^n\E\big|g_i(X)-\E[g_i(X)|X_i]\big|^2.
\end{align*} 
\end{proposition}
\begin{proof}
Let $A_k=\avg_{|v|=k}A(v)$ denote the left-hand side, with $A(v)$ denoting the expression being averaged. Note that
\[
A_n=A([n]) = \sum_{i=1}^n\E\big|g_i(X)-\E[g_i(X)|X_i]\big|^2 = \sum_{i=1}^n\E\Var(g_i(X)\,|\,X_i).
\]
We first prove the $k=2$ case. Apply \eqref{ineq:Manuel} (with $\lambda=b/a$) to the conditional law given $X_i$:
\begin{equation*}
\sum_{j \neq i} \E\Big[\big|\E[g_i(X)\,|\,X_i,X_j] - \E[g_i(X)\,|\,X_i]\big|^2 \,\Big|\,X_i\Big] 	\le \frac{b}{a}\Var\big(g_i(X) \,\big|\,X_i\big).
\end{equation*}
From this we deduce
\begin{align*}
A_2 &= \frac{2}{n(n-1)}\sum_{i=1}^n\sum_{j \neq i} \E\Big| \E[g_i(X)\,|\,X_i,X_j] - \E[g_i(X)\,|\,X_i]\Big|^2 \\
	&\le \frac{2b}{an(n-1)}\sum_{i=1}^n \E\Var\big(g_i(X) \,\big|\,X_i\big) = \frac{2b}{an(n-1)}A_n.
\end{align*}
This proves the $k=2$ case, for every $n\ge 2$.

Let us now turn to the general case. We proceed by induction, supposing that the theorem is true for averages over sets of cardinality $k-1$ for all $(n-1)$-dimensional distributions satisfying the assumptions of Theorem \ref{th:independent}.
For $i \in v \subset [n]$ with $|v|=k$, we compute
\begin{align*}
A(v) &= \sum_{j \in v}\E\big|\E[g_j(X)\,|\,X_v]-\E[g_j(X)\,|\,X_j]\big|^2 \\
	&= \E\big|\E[g_i(X)\,|\,X_v]-\E[g_i(X)\,|\,X_i]\big|^2 + \sum_{j \in v\setminus \{i\}}\E\big|\E[g_j(X)\,|\,X_v]-\E[g_j(X)\,|\,X_{ij}]\big|^2  \\
	&\qquad + \sum_{j \in v\setminus \{i\}}\E \big|\E[g_j(X)\,|\,X_{ij}]-\E[g_j(X)\,|\,X_{j}]\big|^2.
\end{align*}
Average over the choice of $i \in v$ on both sides to get
\begin{align*}
A(v) &= \frac{1}{k}A(v) + \frac{1}{k}\sum_{i,j \in v, \, j \neq i}\E\big|\E[g_j(X)\,|\,X_{v}]-\E[g_j(X)\,|\,X_{ij}]\big|^2  \\
	&\qquad + \frac{1}{k}\sum_{i,j \in v, \, j \neq i} \E\big|\E[g_j(X)\,|\,X_{ij}]-\E[g_j(X)\,|\,X_{j}]\big|^2.
\end{align*}
Rearranging,
\begin{align}
A(v) &= \frac{1}{k-1}\sum_{i,j \in v, \, j \neq i}\E\big|\E[g_j(X)\,|\,X_{v}]-\E[g_j(X)\,|\,X_{ij}]\big|^2 \nonumber \\
	&\qquad + \frac{1}{k-1}\sum_{i,j \in v, \, j \neq i} \E \big|\E[g_j(X)\,|\,X_{ij}]-\E[g_j(X)\,|\,X_{j}]\big|^2 \nonumber \\
	&= \frac{1}{k-1}\sum_{i \in v}\E A(v\setminus \{i\}\,|\,i) + \frac{1}{2(k-1)}\sum_{i,j \in v, \, j \neq i} A(ij), \label{pf:mainindep01}
\end{align}
where we write, for any $u \subset [n] \setminus \{i\}$,
\[
A(u\,|\,i) = \sum_{j \in u}\E\big[\big|\E[g_j(X)\,|\,X_{u \cup \{i\}}]-\E[g_j(X)\,|\,X_{ij}]\big|^2 \,|\,X_i\big].
\]
This is identical to the definition of $A(v \setminus \{i\})$, just under the conditional law of $(X_j)_{j \neq i}$ given $X_i$.
We next average over sets $v$ of size $k$.
Note for each distinct $i,j \in [n]$ that a random subset of $[n]$ of cardinality $k$ contains $i$ with probability $k/n$, and contains both $i$ and $j$ with probability $k(k-1)/n(n-1)$. This implies the elementary identity
\begin{equation}
\avg_{|v|=k}\sum_{i,j \in v, \, j \neq i}a_{ij} = \frac{k(k-1)}{n(n-1)}\sum_{i,j \in [n], \, j \neq i}a_{ij}, \label{eq:averaging-elementary}
\end{equation}
for any elements $a_{ij}$ of a vector space.
Moreover, a random set $v \subset [n]$ of size $k$ conditioned to contain $i$ has the same distribution as $u \cup \{i\}$ where $u$ is a random subset of $[n] \setminus \{i\}$ of size $k-1$.
Using these observations, averaging \eqref{pf:mainindep01} over $|v|=k$ yields
\begin{align}
A_k &\le \frac{1}{k-1}\E\sum_{i=1}^n\avg_{|v|=k}\big(1_{i \in v}  A(v\setminus \{i\}\,|\,i)\big) + \frac{k}{2n(n-1)}\sum_{i,j \in [n], \, j \neq i} A(ij) \nonumber \\
	&= \frac{k}{n(k-1)}\E \sum_{i=1}^n\avg_{u\subset [n]\setminus \{i\},\, |u|=k-1} A(u\,|\,i) + \frac{k}{2}A_2. \label{pf:entropyfactorization1}
\end{align}
Now, apply the induction hypothesis to the conditional measure given $X_i$ to get
\begin{equation}
\avg_{u\subset [n]\setminus \{i\},\, |u|=k-1} A(u\,|\,i) \le \frac{b}{a}\frac{(k-1)(k-2)}{(n-1)(n-2)}\sum_{j \in [n] \setminus \{i\}}\E\big[\big|g_j(X) - \E[g_j(X)\,|\,X_{ij}]\big|^2\,|\,X_i\big]. \label{pf:entropyfactorization2}
\end{equation}
The final term simplifies in expectation when summed over $i$:
\begin{align*}
\E&\sum_{i,j \in [n], \, j \neq i}\E\big[\big|g_j(X) - \E[g_j(X)\,|\,X_{ij}]\big|^2\,|\,X_i\big] \\
	&=  \sum_{i,j \in [n], \, j \neq i}\Big(\E\big| g_j(X) - \E[g_j(X)\,|\,X_{j}] \big|^2 - \E\big| \E[g_j(X)\,|\,X_{ij}] - \E[g_j(X)\,|\,X_{j}] \big|^2 \Big) \\
	&= (n-1)A_n - \sum_{1 \le i < j \le n}A(ij) \\
	&= (n-1)A_n - \frac12 n(n-1)A_2.
\end{align*}
Plugging this into \eqref{pf:entropyfactorization2} and then \eqref{pf:entropyfactorization1}, we get
\begin{align*}
A_k &\le \frac{k}{n(k-1)} \frac{b}{a}\frac{(k-1)(k-2)}{(n-1)(n-2)} \Big((n-1)A_n - \frac12 n(n-1)A_2\Big) +  \frac{k}{2}A_2 \\
	&= \frac{b}{a} \frac{k(k-2)}{n(n-2)}A_n + \frac{k}{2}\Big(1 - \frac{b}{a}\frac{k-2}{n-2} \Big)A_2.
\end{align*}
If the parenthetical term in front of $A_2$ is negative, we discard it and use $(k-2)/(n-2)\le(k-1)/(n-1)$ to complete the proof.
Otherwise, we use the $k=2$ case and $b/a \ge 1$ to get
\begin{align*}
A_k &\le \frac{b}{a} \frac{k(k-2)}{n(n-2)}A_n +  \frac{b}{a}\frac{k}{n(n-1)}\Big(1 - \frac{k-2}{n-2} \Big)A_n = \frac{b}{a}\frac{k(k-1)}{n(n-1)}A_n.
\end{align*}
This completes the proof.

\end{proof}

We will apply Proposition \ref{pr:keyvariancebound} to control the relative Fisher information between $\rho_v$ and the product of its marginals.
The relative Fisher information is defined for two probability measures on a Euclidean space by
\begin{equation}
I(\nu\,|\,\mu) := 4\int \big|\nabla \sqrt{ d\nu/d\mu}\big|^2\,d\mu = \int \big|\nabla \log d\nu/d\mu\big|^2\,d\nu, \label{def:Fisher}
\end{equation}
whenever $\nu \ll \mu$ and $\sqrt{ d\nu/d\mu}$ is weakly differentiable.
Define the ``Fisher total correlation'' analogously to $\HC$ by $\FC(\rho_v) := I(\rho_v\,|\,\otimes_{i \in v}\rho_i)$.

\begin{corollary} \label{co:FisherInfo}
Let $X\sim \rho$ with $\rho$ as in Theorem \ref{th:independent}.
For $2 \le k \le n$,
\[
\avg_{|v|=k} \FC(\rho_v) \le \frac{b}{a} \frac{k(k-1)}{n(n-1)}\FC(\rho).
\]
\end{corollary}
\begin{proof}
Using the marginalization identity $\partial_i\log\rho_v(X_v)=\E[\partial_i\log\rho(X)|X_v]$ for $i \in v \subset [n]$, 
\begin{align*}
\FC(\rho_v) &= \sum_{i \in v} \E \big| \E[\partial_i \log\rho(X)\,|\,X_v] - \E[\partial_i \log\rho(X)\,|\,X_i]\big|^2.
\end{align*}
The claim follows by applying Proposition \ref{pr:keyvariancebound} with $g_i=\partial_i\log \rho$.
\end{proof}

In fact, all of the arguments so far in this section remain valid beyond the log-concave class. Specifically, Corollary \ref{co:FisherInfo} remains valid under just the assumption that $\rho$ and all of its conditional measures satisfy the inequality \eqref{ineq:variance-main} with a uniform constant, i.e.,
\[
\Var\Big(\sum_{i=1}^n f_i(X_i) \,\big|\,X_u \Big) \le \lambda\sum_{i=1}^n\Var(f_i(X_i)\,|\,X_u), \quad a.s., \ \forall u \subset [n].
\]
However, transferring the Fisher information bound to the entropy bound of Theorem \ref{th:independent} requires more. We integrate the bound of Corollary \ref{co:FisherInfo} using de Bruijn's identity, using crucially the fact that the ``condition number'' $b/a$ is nonincreasing along heat flow. First, we give one last somewhat standard lemma, stating that this condition number is inherited by marginals:

\begin{lemma} \label{le:marginalbounds}
Suppose $\rho$ satisfies the assumptions of Theorem \ref{th:independent}. Then, for any $v \subset [n]$, we have the same bounds for the marginals:
\[
aI \le -\nabla^2\log\rho_v(x_v), \qquad -\partial_{ii}\log\rho_v(x_v) \le b, \quad \forall x_v \in \R^v, \ i \in v.
\]
\end{lemma}
\begin{proof}
The lower bound simply says that the marginal of a strongly log-concave measure is again strongly log-concave with the same parameter, and this is a standard consequence of the Pr\'ekopa-Leindler inequality. The upper bound follows from the representation
\[
-\partial_{ii}\log\rho_v(X_v) = -\E[\partial_{ii}\log\rho(X)\,|\,X_v] - \Var(\partial_i\log\rho(X)\,|\,X_v),
\]
which follows by differentiating twice in the marginal formula $\rho_v(x_v)=\int \rho(x) \,dx_{[n] \setminus v}$.
\end{proof}

\begin{proof}[Proof of Theorem \ref{th:independent}]
Note for any probability measure $\mu$ on Euclidean space that $\HC(\mu)=\H(\otimes_i \mu_i)-\H(\mu)$, where $\H(\mu)=-\int\mu\log\mu$ denotes the differential entropy (as long as $\mu$ admits a density with $\mu\log\mu$ being integrable). Let $\mu^t$ denote the convolution of $\mu$ with the centered Gaussian of the same dimension as $\mu$ with covariance matrix $tI$; that is, $\mu^t$ is the solution of the heat equation started from $\mu$. Note that this convolution commutes with marginalization. Let $\I(\mu) = \int|\nabla\log\mu|^2\,d\mu$ denote the Fisher information. The classical de Bruijn identity gives, for $t>0$,
\[
\frac{d}{dt}\H(\mu^t) = \frac12 \I(\mu^t), \qquad \frac{d}{dt}\H(\mu_i^t) = \frac12 \I(\mu_i^t), \ \ i \in [n].
\]
This is valid for measures with finite second moment.
From this we deduce
\begin{equation}
\frac{d}{dt}\HC(\rho_v^t) = - \frac12 \FC(\rho_v^t). \label{pf:deBruijn}
\end{equation}
Before we put this to use, we first prepare some estimates on the log-Hessian of $\rho^t_v$.
A standard calculation for Gaussian convolutions shows that
\begin{equation}
-\nabla^2 \log \rho_v^t(z_v) = \frac{1}{t}I -\frac{1}{t^2}\Cov(X_v\,|\,X_v+\sqrt{t}Z_v=z_v), \label{eq:convolutionHessian}
\end{equation}
where $Z$ is standard Gaussian in $\R^n$, independent of $X \sim \rho$.
The conditional law of $X_v$ given $X_v+\sqrt{t}Z_v=z_v$ has density 
\[
\mu_{z_v}(x_v):=(2\pi t)^{-|v|/2}\frac{\rho_v(x_v)}{\rho^t_v(z_v)}\exp(-|x_v-z_v|^2/2t).
\]
Because $\rho_v$ satisfies the bounds of Lemma \ref{le:marginalbounds}, this conditional density satisfies
\[
\Big(a+\frac{1}{t}\Big)I \le -\nabla^2\log\mu_{z_v}  , \qquad \max_{i \in v} -\partial_{ii} \log\mu_{z_v} \le b+\frac{1}{t}.
\]
The lower bound and the Brascamp-Lieb inequality yield
\[
\Cov(X_v\,|\,X_v+\sqrt{t}Z_v=z_v) \le \Big(a+\frac{1}{t}\Big)^{-1}I.
\]
The upper bound is inherited by each one-dimensional marginal of $\mu_{z_v}$, and then the Cram\'er-Rao inequality implies
\[
\Var(X_i\,|\,X_v+\sqrt{t}Z_v=z_v) \ge \Big(b+\frac{1}{t}\Big)^{-1}, \quad \forall i \in v.
\]
Returning to \eqref{eq:convolutionHessian}, we deduce
\begin{equation*}
\frac{a}{1+at}I \le -\nabla^2 \log \rho_v^t, \qquad -\partial_{ii}\log \rho_v^t \le \frac{b}{1+bt}, \quad \forall i \in v. 
\end{equation*}
We will argue below   that $\HC(\rho_v^t) \to 0$ as $t\to\infty$ and that $\HC(\rho_v^t) \to \HC(\rho_v)$ as $t\downarrow 0$, and so integrating \eqref{pf:deBruijn} yields
\[
\HC(\rho_v) = \frac12\int_0^\infty \FC(\rho_v^t)\,dt.
\]
Averaging over $v$ and applying Corollary \ref{co:FisherInfo} to $\rho^t$ yields
\begin{align*}
\avg_{|v|=k} \HC(\rho_v) &\le \frac12\frac{k(k-1)}{n(n-1)} \int_0^\infty  \frac{b(1+at)}{a(1+bt)} \FC(\rho^t)\,dt \\
	&\le \frac12 \frac{b}{a}\frac{k(k-1)}{n(n-1)} \int_0^\infty   \FC(\rho^t)\,dt \\
	&= \frac{b}{a}\frac{k(k-1)}{n(n-1)} \HC(\rho). 
\end{align*}
It remains to justify the claim that $\HC(\rho_v^t)$ converges to $0$ as $t\to\infty$ and to $\HC(\rho_v)$ as $t\downarrow 0$.
The short-time limit follows from the lower semicontinuity of relative entropy and the inequality $\HC(\rho_v^t) \le \HC(\rho_v)$ for all $t>0$, which is a consequence of  the data processing inequality.
For the long-time limit, use the Pythagorean identity noted in \eqref{eq:pythagoras} to get
\[
\HC(\rho^t_v) \le H(\rho^t_v\,|\,N_k(m_v,tI)),
\]
where $m_v$ is the mean of $\rho_v$, and $N_k(x,tI)$ denotes the $k$-dimensional Gaussian measure with mean $x$ and covariance matrix $tI$. Noting that $\rho^t_v=\int N_k(x_v,tI)\,\rho_v(dx_v)$, we deduce from convexity of entropy that
\[
\HC(\rho^t_v) \le \int H(N_k(x_v,tI)\,|\,N_k(m_v,tI))\,\rho_v(dx_v) = \frac{1}{2t}\E|X_v-m_v|^2 \to 0. \qedhere
\]
\end{proof}

\subsection{Proof of Corollary \ref{co:main-stronger}}
Letting $\H$ again denote the differential entropy,
\begin{equation}
\HC(\rho) = \sum_{i=1}^n \H(\rho_i) - \H(\rho).\label{pf:co1}
\end{equation}
We will bound each term separately.
For the first term, note that the Poincar\'e inequality implies that the variance of $\rho_i$ is at most $1/a$. Thus, 
\begin{equation}
\H(\rho_i) \le \H(N(0,1/a)) = \frac12\log\frac{2\pi e}{a}. \label{pf:co:h(rho_i)}
\end{equation}

Let $\I$ again denote the Fisher information. The entropy power 
obeys an uncertainty principle known as Stam's entropy-Fisher information inequality,
\[
(2\pi e)^{-1}e^{2\H(\rho)/n}\I(\rho) \ge n.
\]
See \cite[Proposition 3.4]{RaginskySason} for proof; in fact, it was discovered by Carlen \cite{carlen1991superadditivity} that this inequality is equivalent to Gross's log-Sobolev inequality for the Gaussian measure.
Using integration by parts followed by the upper bound $-\partial_{ii}\log\rho \le b$ for all $i$, we have $\I(\rho) \le nb$.
Thus, 
\begin{equation*}
\H(\rho) \ge \frac{n}{2}\log\frac{2\pi e n}{\I(\rho)} \ge \frac{n}{2}\log\frac{2\pi e}{b}.
\end{equation*}
Plug this and \eqref{pf:co:h(rho_i)} into \eqref{pf:co1} to get $\HC(\rho) \le (n/2)\log(b/a)$.
 
\hfill\qedsymbol

\subsection{Proof of Corollary \ref{co:LSI+Poincare}}

The proof is a combination of the log-Sobolev inequality for the product measure $\otimes_{i\in[n]}\rho_i$ with the Poincar\'e inequality for the conditional measures given one coordinate. This idea appeared also in \cite{LackerMukherjeeYeung} and \cite{BaptistaMorrisonZahmMarzouk}.
Note first that $-\nabla^2\log\otimes_{i \in [n]}\rho_i \ge aI$, and so Theorem \ref{th:independent} combined with the Bakry-\'Emery log-Sobolev inequality yield
\[
\avg_{|v|=k}\HC(\rho_v) \le \frac{b}{2a^2}\frac{k(k-1)}{n(n-1)}I(\rho\,|\,\otimes_{i \in [n]}\rho_i).
\]
Using $(-\log\rho_i)'(X_i)=\E[-\partial_i\log\rho(X)\,|\,X_i]$ for $X \sim \rho$, we can write
\[
I(\rho\,|\,\otimes_{i \in [n]}\rho_i) = \sum_{i=1}^n \E \, \Var(\partial_i\log\rho(X)\,|\,X_i).
\]
The conditional law of $(X_j)_{j \neq i}$ given $X_i$ is strongly log-concave with the same parameter $a$, and the Poincar\'e inequality thus yields
\[
\Var(\partial_i\log\rho(X)\,|\,X_i) \le \frac{1}{a}\sum_{j \neq i} \E\big[ |\partial_{ij}\log\rho(X)|^2\,|\,X_i\big], \quad a.s.
\]
Combining the above three displays completes the proof. \hfill\qedsymbol

\subsection{Proof of Proposition \ref{pr:GaussianCase}}
This could be deduced quickly from \cite{DevroyeMehrabianReddad}, but we give a quick and more self-contained proof.
Recall that the left-hand side of the claimed inequality is unchanged if $\rho$ is replaced by its pushforward by any map of the form $x\mapsto x_0+Dx$, where $D$ is an invertible diagonal matrix and $x_0 \in \R^n$. We  may thus assume without loss of generality that $\rho$ is centered and that its covariance matrix $\Sigma$ satisfies $\Sigma_{ii}=1$ for all $i$. Then $\Sigma=\Corr(\rho) \le \lambda I$ by assumption, and $\otimes_{i \in v}\rho_i$ is the standard Gaussian in dimension $|v|$.

Letting $\Sigma_v=\Corr(\rho_v)$, we first compute (recalling the elementary identity \eqref{eq:averaging-elementary})
\begin{align*}
\avg_{|v|=k}\tr((\Sigma_v-I)^2) = \avg_{|v|=k}\sum_{i,j \in v, \, i\neq j}\Sigma_{ij}^2 &= \frac{k(k-1)}{n(n-1)}\sum_{i,j \in [n], \, i \neq j}\Sigma_{ij}^2 \\
	&= \frac{k(k-1)}{n(n-1)}\big(\tr \,\Sigma^2 - n\big) \\
		&\le (\lambda-1)\frac{k(k-1)}{n-1} =: D,
\end{align*}
where the last step used $\tr\, \Sigma^2\le\lambda\,\tr\,\Sigma = \lambda n$.
Now, let $S$ be the set of $v \subset [n]$ of size $k$ such that $\Sigma_v \ge I/2$. For any $v\notin S$ of size $k$, we have 
\[
\frac14 < \|\Sigma_v-I\|_{\mathrm{op}}^2 \le \tr( (\Sigma_v-I)^2).
\]
Combining the two preceding inequalities, we deduce  from Markov's inequality that a randomly chosen set $v \subset [n]$ of size $k$ belongs to $S$ with probability at least $1-4 D$.
On the other hand, for any $v \in S$, we estimate the relative Fisher information (defined in \eqref{def:Fisher}) by
\begin{align*}
I(\rho_v\,|\,\otimes_{i \in v}\rho_i) &= \tr\big((\Sigma_v^{-1}-I)^2\Sigma_v\big) \\
&= \tr\big((\Sigma_v-I)\Sigma_v^{-1}(\Sigma_v-I)\big) \\
&\le 2\tr\big((\Sigma_v-I)^2\big).
\end{align*}
By the Pinsker and log-Sobolev inequalities (recall that $\rho_i$ is standard Gaussian), we deduce
\[
\|\rho_v-\otimes_{i \in v}\rho_i\|_{\mathrm{TV}}^2 \le 2\tr\big((\Sigma_v-I)^2\big), \quad v \in S.
\]
Combining the two cases, and using the trivial bound $\|\cdot\|_{\mathrm{TV}} \le 2$, we deduce
\begin{align*}
\avg_{|v|=k}\|\rho_v-\otimes_{i \in v}\rho_i\|_{\mathrm{TV}}^2 &\le 16D + 2\avg_{|v|=k}\Big(1_{v \in S} \tr\big((\Sigma_v-I)^2\big) \Big) \le 18D. 
\end{align*}
This completes the proof. \hfill\qedsymbol

\section{Proof of Theorem \ref{th:main}} \label{se:proof:main}

We give the proof of the TV bound, which is the more involved. In the short Section \ref{se:proof:main-BL} we explain the adaptations for the BL bound.

It suffices to prove the theorem under the additional assumption that $\rho$ has a smooth strictly positive density; if it does not, we convolve with a centered Gaussian with covariance matrix $\epsilon I$, which preserves the upper bound on the correlation matrix, and we can ultimately take the limit $\epsilon \downarrow 0$.

We may assume without loss of generality that $X \sim \rho$ satisfies $\Var(X_i)=1$ for all $i$. Indeed, if this is not the case, then we could replace $X$ with $\widetilde{X}=D^{-1/2}X$ where $D$ is the diagonal matrix formed by $(\Var(X_i))_{i \in [n]}$. The law $\widetilde{\rho}$ of $\widetilde{X}$ then satisfies $\|\rho_v-\otimes_{i \in v}\rho_i\|_{\mathrm{TV}}=\|\widetilde{\rho}_v-\otimes_{i \in v}\widetilde{\rho}_i\|_{\mathrm{TV}}$ for all $v$, so this diagonal transformation does not change the quantity that we are trying to bound. Note that $\Cov(\widetilde{\rho})=\Corr(\rho) \le \lambda I$ by assumption.

With the previous two paragraphs in mind, we assume henceforth that $X\sim\rho$ satisfies the assumptions of Theorem \ref{th:main} and also has a smooth strictly positive density and $\Cov(\rho)=\Corr(\rho)$; in particular, $\Var(X_i)=1$ for all $i$.
We will modify $\rho$ in two ways, first to create an upper bound on $-\nabla^2\log\rho$, and second to create a lower bound.

\subsection{Heat flow}

In this section we regularize $\rho$ by heat flow to create an upper bound on $-\nabla^2\log\rho$.
Let $\rho^s$ be the convolution of $\rho$ with the centered Gaussian of covariance matrix $sI$, and similarly for $\rho^s_v$. Note that Gaussian convolution and marginalization commute.  
Gaussian convolution also preserves log-concavity, as a standard consequence of the Pr\'ekopa-Leindler inequality.
A standard calculation shows that heat flow produces an upper Hessian bound,
\begin{align*}
-\nabla^2\log \rho^s_v (x_v) = \frac{1}{s}I - \frac{1}{s^2}\Cov(X_v\,|\,X_v+\sqrt{s}Z_v=x_v) \le \frac{1}{s}I,
\end{align*}
where $X\sim \rho$ and $Z$ is standard Gaussian in dimension $n$ independent of $X$.
Moreover, we have $\Cov(\rho^s) = \Cov(\rho) + sI$, which implies
\begin{equation}
sI \le \Cov(\rho^s) \le (\lambda +s)I, \qquad \Cov(\rho^s)_{ii} = 1 + s, \ \forall i \in [n]. \label{ineq:heatflow-Cov}
\end{equation}
We need an anisotropic variant of an estimate of Eldan-Klartag \cite[Proposition 10]{EldanKlartag2008}.

\begin{lemma} \label{le:EldanKlartag}
Let $m > 0$.
Let $\mu$ be a smooth strictly positive log-concave probability density function on $\R^k$, and suppose all eigenvalues of $\Cov(\mu)$  are at least $m$. Then
\[
\|\mu-\mu^s\|_{\mathrm{TV}} \lesssim k\sqrt{s/m}.
\]
\end{lemma}
\begin{proof}
We may assume without loss of generality that $\mu$ has mean zero.
We start by noting that
\begin{align*}
\|\mu-\mu^s\|_{\mathrm{TV}} &= \int_{\R^k} \big|\E\, \mu(x + \sqrt{s}Z)-\mu(x)\big|\,dx \\
	&\le \sqrt{s}\int_0^1\E\int_{\R^k} \big|\langle Z, \nabla \mu(x + u\sqrt{s}Z)\rangle\big|\,dx\,du \\
	&= \sqrt{s} \,\E\int_{\R^k} \big|\langle Z, \nabla \mu(x)\rangle\big| \,dx \\
	&= \sqrt{\frac{2s}{\pi}}\int_{\R^k}|\nabla\mu(x)|\,dx.
\end{align*}
This is a simple instance of a much more general inequality on regularization by Markov semigroups due to Ledoux \cite[Inequality (5.5)]{ledoux2004spectral}.
Next, let $A=\Cov(\mu)^{1/2}$ and set $\nu(x)=\det A\mu(Ax)$, so that $\nu$ is the law of $A^{-1}X$ for $X \sim \mu$. Then $\nu$ is isotropic, so we may apply  \cite[Lemma 11]{EldanKlartag2008}, which is where log-concavity is required. It yields
\[
\int |\nabla\nu(x)|\,dx \lesssim  k.
\]
Using $\nabla\mu(x)=A^{-1}\nabla\nu(A^{-1}x)/\det  A$ and $\|A^{-1}\|_{\mathrm{op}} \le 1/\sqrt{m}$, the claim follows.
\end{proof}

Next, we apply this to the marginals of our given measure $\rho$.

\begin{lemma} \label{le:heatflow}
For $1 \le k \le n$ and $s>0$,
\begin{align}
\avg_{|v|=k}\|\rho_v-\rho_v^s\|_{\mathrm{TV}} &\lesssim k\sqrt{s}+(\lambda-1)\frac{k^2}{n},  \label{ineq:heatflow1} \\
\avg_{|v|=k}\|\otimes_{i \in v}\rho_i - \otimes_{i \in v}\rho_i^s\|_{\mathrm{TV}} &\lesssim  k\sqrt{s}. \label{ineq:heatflow2}
\end{align}
\end{lemma}
\begin{proof}
The second inequality is easier. Begin with the standard subadditivity inequality for the total variation distance between product measures:
\[
\|\otimes_{i \in v}\rho_i - \otimes_{i \in v}\rho_i^s\|_{\mathrm{TV}} \le \sum_{i \in v}\|\rho_i -\rho_i^s\|_{\mathrm{TV}}.
\]
Recalling that the one-dimensional marginal $\rho_i$ has variance 1 for each $i$, we may apply Lemma \ref{le:EldanKlartag} with $k=1$ and $m=1$ to  deduce \eqref{ineq:heatflow2}.

The proof of \eqref{ineq:heatflow1} requires more care because we do not have an a priori lower bound on $\Sigma$. 
Recalling that $\Sigma:=\Cov(\rho)$ has all ones on the diagonal and $\Sigma \le \lambda I$, we have
\[
\sum_{i,j \in [n], \, i\neq j} \Sigma_{ij}^2 = \tr(\Sigma^2)-n \le \lambda\, \tr\,\Sigma - n = (\lambda-1)n.
\] 
Thus, recalling the elementary averaging identity \eqref{eq:averaging-elementary},
\begin{align*}
\avg_{|v|=k}\tr( (\Cov(\rho_v)-I)^2) &= \avg_{|v|=k} \sum_{i,j \in v, \, i \neq j}\Sigma_{ij}^2  \\
	&= \frac{k(k-1)}{n(n-1)}\sum_{i,j \in [n], \, i\neq j} \Sigma_{ij}^2 \\
	&\le (\lambda-1)\frac{k(k-1)}{n-1}.
\end{align*}
Let $S$ be the set of subsets $v\subset [n]$ of size $k$ for which $\Cov(\rho_v) \ge I/2$. For any $v\notin S$ of size $k$, we have 
\[
\frac14 < \|\Cov(\rho_v)-I\|_{\mathrm{op}}^2 \le \tr( (\Cov(\rho_v)-I)^2),
\]
Combining the two preceding inequalities, we deduce  from Markov's inequality that a randomly chosen set $v \subset [n]$ of size $k$ belongs to $S$ with probability at least $1-4 (\lambda-1)\frac{k(k-1)}{n-1}$.
Hence, bounding the total variation trivially by 2 on this bad set,
\begin{align*}
\avg_{|v|=k}\|\rho_v-\rho_v^s\|_{\mathrm{TV}} &\le \avg_{|v|=k}\big(1_{v \in S}\|\rho_v-\rho_v^s\|_{\mathrm{TV}}\big) + 8 (\lambda-1)\frac{k(k-1)}{n-1}.
\end{align*}
For $v \in S$, we have $\Cov(\rho_v) \ge I/2$ by definition, and by Lemma \ref{le:EldanKlartag} the first term is $\lesssim k\sqrt{2s}$. This completes the proof.

\end{proof}

\subsection{Stochastic localization} \label{se:stochasticlocalization}
In this section we regularize $\rho^s$ to make it \emph{strongly} log-concave. There are multiple ways to do this.
The simplest is to multiply by a Gaussian, replacing $\rho^s$ with the probability measure proportional to $e^{-t|x|^2/2}\rho^s(dx)$, for $t>0$.
We instead follow a more complicated but apparently more efficient construction, known as stochastic localization. Both approaches seem feasible for proving versions of Theorem \ref{th:main}, but we were able to achieve better dependence on $k$ and $n$ using stochastic localization.

Fix $s>0$ for now, and let $C=\Cov(\rho^s)^{-1/2}$.
Let $(\rho^{s,t})_{t \ge 0}$ be the stochastic localization with a modified covariation,  given as the solution of the family of stochastic differential equations
\[
d\rho^{s,t}(x) = \rho^{s,t}(x)\langle C(x-a_t), dB_t\rangle,  \quad x \in \R^n, \quad a_t = \int y\,\rho^{s,t}(y)\,dy, 
\]
initialized at $\rho^{s,0}=\rho^s$. Here $B$ is an $n$-dimensional Brownian motion, and we refer to \cite[Section 2]{eldan2013thin} for a rigorous construction of this process.
A standard calculation, applying It\^o's formula to $d\log \rho^{s,t}(x)$, leads to
\[
\rho^{s,t}(x) = \rho^s(x)\exp\bigg(\int_0^t \langle C(x-a_r), dB_r\rangle - \frac12 \int_0^t |C(x-a_r)|^2\,dr\bigg).
\]
This formula implies
\[
-\nabla^2\log \rho^{s,t} = -\nabla^2\log \rho^s + t\Cov(\rho^s)^{-1}.
\]
Using \eqref{ineq:heatflow-Cov} and $0 \le -\nabla^2\log \rho^s \le (1/s)I$,  we deduce that
\[
\frac{t}{\lambda+s}I \le  -\nabla^2\log \rho^{s,t} \le \frac{1+t}{s} I.
\]
The ratio of the upper and lower bounds is
\begin{equation}
\lambda_t :=  \Big(1 + \frac{\lambda}{s}\Big)\Big(1 + \frac{1}{t}\Big). \label{pf:def:lambda_t}
\end{equation}
We next compute some dynamics.
Marginalizing, we have
\begin{align*}
d\rho^{s,t}_v(x_v) &=  \rho^{s,t}_v(x_v) \big\langle\E_{\rho^{s,t}}\big[C(X-a_t)\,|\,X_v=x_v\big] , \, dB_t\big\rangle, \\
d\rho^{s,t}_i(x_i) &= \rho^{s,t}_i(x_i) \big\langle\E_{\rho^{s,t}}\big[C(X-a_t)\,|\,X_i=x_i\big] , \, dB_t\big\rangle.
\end{align*}
Here $\E_{\rho^{s,t}}$ indicates that $X\sim \rho^{s,t}$, with $\rho^{s,t}$ being a random measure, so the resulting quantity remains random; below, we will reserve the subscript-free symbol $\E$ for expectation over this remaining randomness.
By the It\^o product rule, and setting $u_i^t(x_i)=\E_{\rho^{s,t}}[C(X-a_t)\,|\,X_i=x_i]$, we have
\begin{align*}
d\otimes_{i \in v}\rho^{s,t}_i(x_v) &= \sum_{i \in v} \big(\otimes_{\ell \in v}\rho^{s,t}_\ell(x_\ell)\big)\big\langle u^t_{i}(x_i), \, dB_t\big\rangle \\
	&\quad + \sum_{i,j \in v, \, i < j}\big(\otimes_{\ell \in v}\rho^{s,t}_\ell(x_\ell)\big) \big\langle u^t_{i}(x_i),  u^t_{j}(x_j)\big\rangle\,dt .
\end{align*}
Take the difference, integrate against a test function $|f_v| \le 1$, and then apply Cauchy-Schwarz:
\begin{align*}
\Big|\frac{d}{dt}\E\langle  \rho^{s,t}_v-\otimes_{i \in v}\rho^{s,t}_i, f_v\rangle\Big|^2 &= \Big|\E\sum_{i,j \in v, \, i < j}\E_{\otimes_{\ell \in v}\rho^{s,t}_\ell}\big[ f_v(X_v) \big\langle u^t_{i}(X_i),  u^t_{j}(X_j)\big\rangle\big]\Big|^2 \\
	&\le  \E\,\E_{\otimes_{\ell \in v}\rho^{s,t}_\ell}\bigg[  \Big( \sum_{i,j \in v, \, i < j}\big\langle u^t_{i}(X_i), u^t_{j}(X_j)\big\rangle \Big)^2\bigg]  \\
	&=  \E\sum_{i,j \in v, \, i < j} \E_{\otimes_{\ell \in v}\rho^{s,t}_\ell}\big[ \big\langle u^t_{i}(X_i), u^t_{j}(X_j)\big\rangle^2\big]  \\
	&= \E\sum_{i,j \in v, \, i < j} \tr(K^t_iK^t_j) ,
\end{align*}
where we set $K^t_i = \int u^t_i(x_i)u^t_i(x_i)^\top \,\rho^{s,t}_i(dx_i)$, and where the second to last identity used the fact that $\langle u^t_{i}(x_i), u^t_{j}(x_j)\rangle$ and $\langle u^t_{i'}(x_{i'}), u^t_{j'}(x_{j'})\rangle$ are orthogonal in $L^2(\otimes_{\ell \in v}\rho^{s,t}_\ell)$ whenever $\{i,j\} \neq \{i',j'\}$.
Now, averaging over $v$ (recalling the elementary identity \eqref{eq:averaging-elementary}),
\begin{equation}
\Big|\avg_{|v|=k}\frac{d}{dt}\E\langle  \rho^{s,t}_v-\otimes_{i \in v}\rho^{s,t}_i, f_v\rangle\Big|^2 \le  \frac{k(k-1)}{n(n-1)} \E \sum_{1 \le i < j \le n} \tr(K^t_iK^t_j). \label{pf:stochloc3}
\end{equation}
To control the trace, we apply the inequality \eqref{ineq:Manuel} to $g(x)=\langle \theta, C(x-a_t)\rangle$ for a given $\theta \in \mathbb{S}^{n-1}$ to find
\begin{align*}
\sum_{i=1}^n \langle\theta,K^t_i\theta\rangle &= \sum_{i=1}^n\E_{\rho^{s,t}}  \langle \theta,u^t_i(X_i)\rangle^2 \\
	  &= \sum_{i=1}^n\E_{\rho^{s,t}}|\E_{\rho^{s,t}}[\langle \theta, C(X-a_t)\rangle\,|\,X_i]|^2 \\
	&\le \lambda_t \E_{\rho^{s,t}} \langle \theta, C(X-a_t)\rangle^2 \\
	&= \lambda_t \langle\theta, C\Cov(\rho^{s,t})C \theta\rangle.
\end{align*}
In other words, we have in semidefinite order that
\[
\sum_{i=1}^n K^t_i \le \lambda_t C\Cov(\rho^{s,t})C.
\]
Alternatively, Jensen's inequality yields for each $i$ that
\[
K^t_i = \E_{\rho^{s,t}}\big[\E_{\rho^{s,t}}[C(X-a_t)\,|\,X_i]^{\otimes 2}\big] \le \E_{\rho^{s,t}}\big[ (C(X-a_t))^{\otimes 2}\big] = C\Cov(\rho^{s,t})C.
\]
Summing over $i \in [n]$ and combining with the preceding inequality yields
\[
\sum_{i=1}^n K^t_i \le (n \wedge \lambda_t) C\Cov(\rho^{s,t})C.
\]
Returning to \eqref{pf:stochloc3}, and including the diagonal term $i=j$,  we get\footnote{We write $\tr \, A^2$ for the trace of the square of a matrix, versus $(\tr \, A)^2$ for the square of the trace.}
\begin{align*}
\Big|\frac{d}{dt}\avg_{|v|=k}\E\langle  \rho^{s,t}_v-\otimes_{i \in v}\rho^{s,t}_i, f_v\rangle\Big|^2 &\le \frac{k(k-1)}{2n(n-1)} \E\,\tr\Big( \sum_{i=1}^n K^t_i\Big)^2 \\
	&\le \frac{k(k-1)}{2n(n-1)}(n \wedge \lambda_t)^2 \E\,\tr\,( C\Cov(\rho^{s,t})C)^2 .
\end{align*}
We will control the final trace term using Guan's bound, but to do so we must first transform $\rho^{s,t}$ to isotropic position.
If we define $\overline{p}^t(x) = |\det C|^{-1}\rho^{s,t}(C^{-1}x)$, so that $\overline{p}^t$ is the pushforward of $\rho^{s,t}$ by the linear map $C$, then $\overline{a}_t:=\int y\,\overline{p}^t(y)\,dy$ satisfies $\overline{a}_t=Ca_t$, and so $\overline{p}^t$ satisfies
\[
d\overline{p}^t(x)=\overline{p}^t(x) \langle x-\overline{a}_t,dB_t\rangle.
\]
Recalling that $C=\Cov(\rho^s)^{-1/2}$, we deduce that $\overline{p}^0$ is isotropic.  Hence,  Guan's bound \cite[Lemma 2.1]{guan2024note} applies  (see \cite[Lemma 2.3]{klartag2025affirmative} for the extension to arbitrary times), yielding
\[
\E\,\tr( C\Cov(\rho^{s,t})C)^2 = \E \,\tr\,\Cov(\overline{p}^t)^2 \lesssim n, \qquad \forall t > 0.
\]
We have thus finally shown that
\begin{equation*}
\Big|\frac{d}{dt}\avg_{|v|=k}\E\langle  \rho^{s,t}_v-\otimes_{i \in v}\rho^{s,t}_i, f_v\rangle\Big| \lesssim  \frac{k}{\sqrt{n}}(n \wedge \lambda_t).
\end{equation*}
Integrate up to time 1, recalling that $\rho^{s,0}=\rho^s$, to get
\begin{equation*}
\avg_{|v|=k} \langle  \rho^s_v-\otimes_{i \in v}\rho^s_i, f_v\rangle - \avg_{|v|=k}\E\langle  \rho^{s,1}_v-\otimes_{i \in v}\rho^{s,1}_i, f_v\rangle  \lesssim    \frac{k}{\sqrt{n}}\int_0^1 (n \wedge \lambda_t)\,dt.
\end{equation*}
Taking a supremum over $f_v$ and using Pinsker's inequality, we get
\begin{equation*}
\avg_{|v|=k} \| \rho^s_v - \otimes_{i \in v}\rho^s_i\|_{\mathrm{TV}} \lesssim \sqrt{\avg_{|v|=k}\E \HC( \rho^{s,1}_v)} +   \frac{k}{\sqrt{n}}\int_0^1 (n \wedge \lambda_t)\,dt.
\end{equation*}
To control the integral, we split it at $t=1/n$ and recall the definition of $\lambda_t$ from \eqref{pf:def:lambda_t}:
\[
\int_0^1 (n \wedge \lambda_t)\,dt \le 1 + \Big(1 + \frac{\lambda}{s}\Big)\int_{1/n}^1  \Big(1 + \frac{1}{t}\Big)\,dt  \lesssim \Big(1 + \frac{\lambda}{s}\Big)\log n.
\]
Applying this and Corollary \ref{co:main-stronger} yields
\begin{equation*}
    \avg_{|v|=k} \| \rho^s_v - \otimes_{i \in v}\rho^s_i\|_{\mathrm{TV}} \lesssim  \sqrt{\lambda_1\log\lambda_1} \frac{k}{\sqrt{n}} +  \frac{k}{\sqrt{n}}\Big(1 + \frac{\lambda}{s}\Big)\log n.
\end{equation*}
Using $\sqrt{\lambda_1 \log \lambda_1} \le \lambda_1=2(1+\lambda/s)$, this simplifies to
\begin{equation}
    \avg_{|v|=k} \| \rho^s_v - \otimes_{i \in v}\rho^s_i\|_{\mathrm{TV}} \lesssim   \frac{k}{\sqrt{n}}\Big(1 + \frac{\lambda}{s}\Big)\log n.  \label{pf:stochloc1}
\end{equation}

\subsection{Completing the proof of the TV estimate} \label{se:completingtheproof}

Combining \eqref{pf:stochloc1} with Lemma \ref{le:heatflow} and the triangle inequality, we get
\begin{align*}
A_k := \avg_{|v|=k} \|\rho_v - \otimes_{i \in v}\rho_i\|_{\mathrm{TV}} &\le \avg_{|v|=k}\|\rho_v-\rho^s_v\|_{\mathrm{TV}} + \avg_{|v|=k} \| \rho^s_v-\otimes_{i \in v}\rho^s_i\|_{\mathrm{TV}} \\
	&\quad + \avg_{|v|=k}\|\otimes_{i \in v}\rho^s_i - \otimes_{i \in v}\rho_i\|_{\mathrm{TV}} \\
	&\lesssim k\sqrt{s} + \lambda\frac{k^2}{n} + \frac{k}{\sqrt{n}}\Big(1+\frac{\lambda}{s}\Big) \log n.
\end{align*}
Taking $s=(\lambda  \log n/\sqrt{n})^{2/3}$ yields
\[
A_k \lesssim \frac{\lambda^{1/3} k  (\log n)^{1/3}}{n^{1/6}} + \frac{k }{\sqrt{n}}\log n + \lambda\frac{k^2}{n}.
\]

Letting $R=\lambda^{1/3} k  (\log n)^{1/3}/n^{1/6}$ denote the first term, we can write the second and third as
\[
\frac{k }{\sqrt{n}}\log n = \frac{R^3}{\lambda k^2}, \qquad \lambda\frac{k^2}{n} = \frac{R^3}{k \sqrt{n}\log n}.
\]
A correlation matrix has diagonal entries all equal to 1, so $1 \le \|\Corr(\rho)\|_{\mathrm{op}} \le \lambda$. Hence, these two terms are $\lesssim R^3$. If $R \le 1$, then $R^3 \le R$ and $A_k \lesssim R$ proves the claimed TV  estimate of Theorem \ref{th:main}. If $R \ge 1$, then  the trivial bound $A_k \le 2 \lesssim R$ applies.

\subsection{Proof of Theorem \ref{th:main}, BL estimate} \label{se:proof:main-BL}

The proof of the BL bound is a simple adaptation of the argument above. Note that we do not start by applying a diagonal transformation as we did in the TV case; we simply work with the assumption $\Cov(\rho) \le \lambda I$. In place of Lemma \ref{le:heatflow}, we need only the easy inequalities
\begin{align}
\avg_{|v|=k}\BL(\rho_v,\rho_v^s) \le \sqrt{sk}, \qquad \avg_{|v|=k}\BL(\otimes_{i \in v}\rho_i , \otimes_{i \in v}\rho_i^s) \le\sqrt{sk}. \label{ineq:heatflow-BL}
\end{align} 
Indeed, both follow from a simple coupling argument: If $Y$ is a $k$-dimensional random vector and $Z$ is an independent standard normal of dimension $k$, then for  1-Lipschitz $f$ we have
\[
\E f(Y) - \E f(Y+\sqrt{s}Z) \le \sqrt{s}\E|Z| \le \sqrt{sk}.
\]
The stochastic localization analysis of Section \ref{se:stochasticlocalization} remains valid for $\rho$ satisfying $\Cov(\rho) \le \lambda I$.
Using BL instead of TV, the triangle inequality combined with \eqref{ineq:heatflow-BL} and \eqref{pf:stochloc1} yields
\begin{equation*}
\widehat{A}_k := \avg_{|v|=k} \BL(\rho_v , \otimes_{i \in v}\rho_i) \lesssim  \sqrt{sk} + \frac{k}{\sqrt{n}}\Big(1+\frac{\lambda}{s}\Big) \log n.
\end{equation*}
Taking $s= (\lambda \sqrt{k}\log n/\sqrt{n})^{2/3}$ yields
\begin{equation}
\widehat{A}_k\lesssim k^{2/3} (\lambda  \log n/\sqrt{n})^{1/3}  + \frac{k}{\sqrt{n}}  \log n. \label{pf:BLbound1}
\end{equation}
Let $R= k^{2/3} (\lambda  \log n/\sqrt{n})^{1/3}$ denote the first term. Then the second can be written as 
\[
\frac{k}{\sqrt{n}}\log n = \frac{R^3}{\lambda k}.
\]
If $R^2 \le \lambda k$, then this is at most $R$, and we deduce $\widehat{A}_k \lesssim  R$, which is the BL estimate of Theorem \ref{th:main}. If instead $R^2 > \lambda k$, then we prove $\widehat{A}_k \lesssim R$ by a cruder argument.
Letting $m$ denote the mean of $\rho$, we have
\[
\BL(\rho_v,\delta_{m_v}) \le \int |x_v-m_v|\,\rho(dx) \le (\tr \,\Cov(\rho_v))^{1/2} \le \sqrt{\lambda k}.
\]
Similarly for the product measure $\otimes_{i \in v}\rho_i$ which has the same mean $m_v$. By the triangle inequality, $\widehat{A}_k \le 2\sqrt{\lambda k} \le 2R$, and the proof is complete. \hfill \qedsymbol

\begin{remark}
The use of Guan's bound in the argument above is not entirely essential.  
For instance, the less sharp estimates of Lee-Vempala \cite{LeeVempala} could lead to a similar final result but with worse exponents on $k$ and $n$. However, the easiest bound $\E\, \tr\,\Cov(\overline{p}^t)^2 \le n/t$ coming from the $t$-log-concavity of $\overline{p}^t$ does not yield anything informative.
\end{remark}

\section{Proof of Theorem \ref{th:BLdistance} } \label{se:proof:BLdistance}

We first prove a variant of Theorem \ref{th:BLdistance} in which the bounded-Lipschitz distance is replaced by a much weaker one.
For two probability measures $\mu,\nu$ on $\R^k$, and $t>0$, define
\begin{align*}
\MMD_t^2(\mu,\nu) := \E\bigg|\int_{\R^k} \exp\Big( \iota  \sum_{j=1}^k x_j Z_j\Big)\,(\mu-\nu)(dx)\bigg|^2,
\end{align*}
where $Z_1,\ldots,Z_k$ are iid from $\gamma_t=N(0,t)$, and $\iota $ is the imaginary unit.
Note that this can be written as
\begin{align*}
\MMD_t^2(\mu,\nu) &= \E \int_{\R^k}\int_{\R^k} \exp\Big( \iota \sum_{j=1}^k (x_j-x'_j) Z_j\Big)\,(\mu-\nu)(dx)(\mu-\nu)(dx') \\
	&=  \int_{\R^k}\int_{\R^k} \exp\Big(-\frac{t}{2}|x-x'|^2\Big)\,(\mu-\nu)(dx)(\mu-\nu)(dx') .
\end{align*}
This is the \emph{maximum mean discrepancy} associated with the Gaussian kernel $(x,y) \mapsto e^{-t|x-y|^2/2}$. It is known to metrize the weak topology. 
When measuring the distance between a measure and the product of its marginals, this distance gives rise to what is known in statistics as the \emph{Hilbert-Schmidt independence criterion} \cite{gretton2005measuring,pfister2018kernel}.

\begin{proposition} \label{pr:MMD}
Every probability measure $\rho$ on $\R^n$ satisfies, for all $2 \le k \le n$ and $t > 0$,
\[
\avg_{|v|=k} \MMD_t(\rho_v,\otimes_{i \in v}\rho_i) \le \sqrt{8t\CPIA(\rho)}\frac{k-1}{\sqrt{n-1}}.
\]
\end{proposition}
\begin{proof}
For sets $v\subset [n]$ define the Fourier transforms of the marginals by
\begin{equation}
\widehat{\rho}_v(z_v) = \int_{\R^n} \exp\Big(\iota \sum_{j \in v} z_j x_j\Big)\,\rho(dx). \label{def:marginalFourier}
\end{equation}
Fix for now a set $v=\{i_1,\ldots,i_k\} \subset [n]$ of cardinality $k$, and let $v_\ell=\{i_1,\ldots,i_{\ell}\}$ for $\ell=1,\ldots,k$.
Applying the triangle inequality to a telescoping sum and noting that $|\widehat{\rho}_j|\le 1$ yields
\begin{equation*}
\MMD_t(\rho_v,\otimes_{i \in v}\rho_i) \!=\! \bigg(\E\Big| \widehat{\rho}_v(Z_v) - \prod_{j \in v}\widehat{\rho}_j(Z_j)\Big|^2\bigg)^{1/2} \!\!\!  \le \sum_{\ell=2}^k \Big(\E\big| \widehat{\rho}_{v_\ell}(Z_{v_\ell}) - \widehat{\rho}_{v_{\ell-1}}(Z_{v_{\ell-1}})\widehat{\rho}_{i_\ell}(Z_{i_\ell})\big|^2\Big)^{1/2} .
\end{equation*}
If $i_1,\ldots,i_\ell$ are chosen uniformly at random without replacement from $[n]$, then the conditional distribution of $i_\ell$ given $v_{\ell-1}$ is uniform over $[n] \setminus v_{\ell-1}$. 
Using the triangle inequality followed by this fact,
\begin{align*}
\avg_{|v|=k}\MMD_t(\rho_v,\otimes_{i \in v}\rho_i) &\le \sum_{\ell=2}^k \avg_{|u|=\ell-1}\bigg( \frac{1}{n-\ell+1}\sum_{ j \in [n] \setminus u} \E\big| \widehat{\rho}_{u \cup \{j\}}(Z_{u \cup \{j\}}) - \widehat{\rho}_u(Z_u)\widehat{\rho}_j(Z_j)\big|^2\bigg)^{1/2}.
\end{align*}
We claim that, for all nonempty $u \subset [n]$, we have
\begin{equation}
\sum_{ j \in [n] \setminus u} \E\big| \widehat{\rho}_{u \cup \{j\}}(Z_{u \cup \{j\}}) - \widehat{\rho}_u(Z_u)\widehat{\rho}_j(Z_j)\big|^2 \le 2t \CPIA(\rho) . \label{pf:MMD1}
\end{equation}
Once this is established, we deduce from the preceding inequality that
\[
\avg_{|v|=k}\MMD_t(\rho_v,\otimes_{i \in v}\rho_i) \le \sqrt{2t\CPIA(\rho)}\sum_{\ell=2}^k  \frac{1}{\sqrt{n-\ell+1}}.
\]
The proof of the proposition is then completed by bounding the summation by
\[
\sum_{\ell=n-k+1}^{n-1}  \frac{1}{\sqrt{\ell}} \le \int_{n-k}^{n-1} \frac{ds}{\sqrt{s}} = 2\frac{k-1}{\sqrt{n-1}+\sqrt{n-k}}.
\]

It remains to prove \eqref{pf:MMD1}.
For each $z_u=(z_\ell)_{\ell \in u}$ and $z \in\R$, define
\[
\varphi^{z_u,z}(x):=\sum_{j \in [n]\setminus u}\varphi^{z_u,z}_j(x_j), \quad \varphi^{z_u,z}_j(x_j):=\overline{\widehat{\rho}_{u \cup \{j\}}(z_u,z)-\widehat{\rho}_u(z_u)\widehat{\rho}_j(z) }\big( e^{\iota z x_j} - \widehat{\rho}_j(z)\big).
\]
Here the notation  $\widehat{\rho}_{u \cup \{j\}}(z_u,z)$ means that $z$ is used for the index-$j$ argument of $\widehat{\rho}_{u \cup \{j\}}$.
\begin{align*}
A_u(z_u,z) := \sum_{ j \in [n] \setminus u} \big| \widehat{\rho}_{u \cup \{j\}}(z_u,z) - \widehat{\rho}_u(z_u)\widehat{\rho}_j(z)\big|^2 &= \int \varphi^{z_u,z}(x)e^{\iota \sum_{\ell \in u} z_\ell x_\ell}\,\rho(dx).
\end{align*}
The point of this definition is that the left-hand side of \eqref{pf:MMD1} is precisely $\E A_u(Z_u,Z)$, for $Z\sim N(0,t)$ independent of $Z_u$.
For each $z_u=(z_\ell)_{\ell \in u}$ and $z \in\R$, note that $\varphi^{z_u,z}$ has $\rho$-mean zero. Hence,  by Cauchy-Schwarz,
\begin{align*}
A_u(z_u,z)^2 &\le  \int |\varphi^{z_u,z}(x)|^2\,\rho(dx) = \Var\big(\varphi^{z_u,z}(X)\big) = \Var_\rho\bigg(\sum_{ j \in [n] \setminus u}\varphi^{z_u,z}_j(X_j) \bigg) \\
	&\le 2\CPIA(\rho) \sum_{j \in [n]\setminus u} \| |(\varphi^{z_u,z}_j)'| \|_{L^\infty(\R)}^2.
\end{align*}
The last step used the definition of $\CPIA(\rho)$ from \eqref{ineq:Poincare-additive}, which was given for real-valued functions, and the additional factor of 2 comes from the easy extension to complex-valued functions (applying \eqref{ineq:Poincare-additive} separately to the real and imaginary part).  To control these derivatives, we compute
\begin{align*}
\big|(\varphi^{z_u,z}_j)'(x_j)\big|^2 % &=  
	&= |z|^2 \big|\widehat{\rho}_{u \cup \{j\}}(z_u,z)-\widehat{\rho}_u(z_u)\widehat{\rho}_j(z)\big|^2.
\end{align*}
Combining this with the preceding inequality yields
\begin{align*}
A_u(z_u,z)^2 &\le  2\CPIA(\rho)|z|^2 \sum_{j \in [n]\setminus u} \big|\widehat{\rho}_{u \cup \{j\}}(z_u,z)-\widehat{\rho}_u(z_u)\widehat{\rho}_j(z)\big|^2 \\
	&= 2\CPIA(\rho)|z|^2 A_u(z_u,z).
\end{align*}
Hence, the left-hand side of \eqref{pf:MMD1} is
\[
\E A_u(Z_u,Z) \le 2\CPIA(\rho)\E|Z|^2 = 2t\CPIA(\rho). \qedhere
\]
\end{proof}

We are now ready to prove Theorem \ref{th:BLdistance}. The main idea is to evolve both $\rho_v$ and $\otimes_{i\in v}\rho_i$ by heat flow. On the one hand, this changes their BL distance by at most $2\sqrt{sk}$. On the other hand, the BL distance between  these smoothed measures can be related to the $\MMD_{1/2s}$ distance between the original measures (up to a truncation error which we will ignore in this paragraph), which contributes a factor of $1/\sqrt{s}$ after applying Proposition \ref{pr:MMD}. This trade-off is optimized in $s$. A similar strategy appeared in the \cite[Section 2.5]{VayerGribonval} which developed methods for estimating Wasserstein distances in terms of kernel norms, uniformly over certain classes of measures.

\begin{proof}[Proof of Theorem \ref{th:BLdistance}]
The finiteness of $\CPIA(\rho)$ implies that $\rho$ has finite second moments; this follows from applying \eqref{ineq:Poincare-additive} to a sequence of bounded functions approximating $f_i(x)=x$.
We may then assume without loss of generality that $\rho$ is centered.
Fix $v \subset [n]$ of size $k$.
Let $s>0$, to be chosen later.
As before, let $\rho^s$ denote the convolution of $\rho$ with the centered Gaussian with covariance matrix $sI$, and similarly for $\rho^s_v$.
An easy coupling argument (as in Section \ref{se:proof:main-BL}) yields
\[
\BL^2(\rho_v,\rho^s_v) \le sk, \qquad \BL^2(\otimes_{i \in v}\rho_i,\otimes_{i \in v}\rho^s_i) \le sk.
\]
Then
\begin{equation}
\BL(\rho_v,\otimes_{i \in v}\rho_i) \le 2\sqrt{sk} + \BL(\rho_v^s,\otimes_{i \in v}\rho_i^s). \label{pf:MMD-to-BL0}
\end{equation}
To estimate the remaining BL distance, pick an arbitrary 1-Lipschitz function $f : \R^k \to \R$ with $|f| \le 1$.
Pick $r>0$, to be specified later.
Then
\begin{equation}
\Big|\int f\,d\big(\rho_v^s-\otimes_{i \in v}\rho_i^s\big)\Big| \le \mathrm{Vol}(B_r)^{1/2}\|\rho_v^s-\otimes_{i \in v}\rho_i^s\|_{L^2(\R^k)} + \int_{B_r^c}|\rho_v^s(x)-\otimes_{i \in v}\rho_i^s(x)|\,dx, \label{pf:MMD-to-BL1}
\end{equation}
where $\mathrm{Vol}$ denotes Lebesgue measure, and $B_r$ is the centered ball of radius $r$ in $\R^k$.
To control the second term, note that the inequality \eqref{ineq:Poincare-additive} implies $\Var(X_i) \le \CPIA(\rho)$ for each $i$. Thus,
\begin{align*}
\int_{\R^k} |x|^2|\rho_v^s(x)-\otimes_{i \in v}\rho_i^s(x)|\,dx &\le 2\sum_{i \in v} \int_{\R} z^2\,\rho_i^s(dz) \\
	&= 2\sum_{i \in v} \int_{\R} z^2\,\rho_i(dz) + 2ks \\
	&\le 2k(\CPIA(\rho) + s).
\end{align*}
By Markov's inequality,
\begin{equation}
\int_{B_r^c}|\rho_v^s(x)-\otimes_{i \in v}\rho_i^s(x)|\,dx \le \frac{2k}{r^2}(\CPIA(\rho) + s).  \label{pf:MMD-to-BL2}
\end{equation}
To estimate the first term in \eqref{pf:MMD-to-BL1} in terms of $\MMD_t$, we naturally pass to Fourier transforms. Indeed, let $\widehat{\rho}_v(z)$ again denote the Fourier transform of $\rho_v$ as in \eqref{def:marginalFourier}, and note that the Fourier transform of $\otimes_{i \in v}\rho_i$ is $z \mapsto \prod_{i \in v}\widehat{\rho}_i(z_i)$. The Fourier transform of $\rho_v^s-\otimes_{i \in v}\rho_i^s$ is then
\[
z\mapsto \Big(\widehat{\rho}_v(z)-\prod_{i \in v}\widehat{\rho}_i(z_i)\Big)e^{-s|z|^2/2}.
\]
Hence, by Plancherel,
\begin{equation*}
\|\rho_v^s-\otimes_{i \in v}\rho_i^s\|_{L^2(\R^k)}^2 = (2\pi)^{-k}\int_{\R^k} \Big|\widehat{\rho}_v(z)-\prod_{i \in v}\widehat{\rho}_i(z_i)\Big|^2e^{-s|z|^2}\,dz.
\end{equation*}
An equivalent formula for $\MMD_t^2(\mu,\nu)$ for probability measures $\mu$ and $\nu$ on $\R^k$ is
\[
\MMD_t^2(\mu,\nu) = (2\pi t)^{-k/2}\int_{\R^k}\big|\widehat{\mu}(z)-\widehat{\nu}(z)\big|^2e^{-|z|^2/2t}\,dz.
\]
If we take $s=1/2t$, comparing the last two formulas yields
\begin{equation*}
\|\rho_v^{1/2t}-\otimes_{i \in v}\rho_i^{1/2t}\|_{L^2(\R^k)}^2 = (t/2\pi)^{k/2}\MMD_t^2(\rho_v,\otimes_{i \in v}\rho_i).
\end{equation*}
This identity is a special case of \cite[Remark 5(f)]{sriperumbudur2010hilbert}.
Using $\mathrm{Vol}(B_r)=\pi^{k/2}r^k/\Gamma(1+(k/2))$ and $\Gamma(1+(k/2)) \ge (k/2e)^{k/2}$, we get
\begin{equation*}
\mathrm{Vol}(B_r)^{1/2}\|\rho_v^{1/2t}-\otimes_{i \in v}\rho_i^{1/2t}\|_{L^2(\R^k)} \le (etr^2/k)^{k/4}\MMD_t(\rho_v,\otimes_{i \in v}\rho_i).
\end{equation*}
Plugging this and \eqref{pf:MMD-to-BL2} into \eqref{pf:MMD-to-BL1}, and then optimizing over $f$ and recalling \eqref{pf:MMD-to-BL0}, we find (for $s=1/2t$)
\begin{equation}
\BL(\rho_v,\otimes_{i \in v}\rho_i) \le \sqrt{\frac{2k}{t}} + \frac{2k}{r^2}\Big(\CPIA(\rho) + \frac{1}{2t}\Big) +  \Big(\frac{etr^2}{k}\Big)^{k/4}\MMD_t(\rho_v,\otimes_{i \in v}\rho_i).
\end{equation}
Average over $v$ of size $k$ and use Proposition \ref{pr:MMD} to get
\begin{equation*}
A_k := \avg_{|v|=k}\BL(\rho_v,\otimes_{i \in v}\rho_i) \le \!\!\sqrt{\frac{2k}{t}} + \frac{2k}{r^2}\Big(\CPIA(\rho) + \frac{1}{2t}\Big) +  \Big(\frac{etr^2}{k}\Big)^{k/4}\!\!\!\!\sqrt{8t\CPIA(\rho)}\frac{k-1}{\sqrt{n-1}}.
\end{equation*}
An elementary computation shows for real numbers $a,b>0$ that
\begin{equation*}
\inf_{r>0} \Big(\frac{a}{r^2} + br^{k/2}\Big)  = a^{\frac{k}{k+4}}b^{\frac{4}{k+4}}(k/4)^{\frac{4}{k+4}}\Big(1 + \frac{4}{k}\Big) \le 5a^{\frac{k}{k+4}}b^{\frac{4}{k+4}}(k/4)^{\frac{4}{k+4}}.
\end{equation*}

This yields
\begin{equation*}
A_k \le \sqrt{\frac{2k}{t}} + 5\Big[2k\Big(\CPIA(\rho) + \frac{1}{2t}\Big)\Big]^{\frac{k}{k+4}} \Big[ (et/k)^{k/4}\sqrt{8t\CPIA(\rho)}\frac{k-1}{\sqrt{n-1}}\Big]^{\frac{4}{k+4}}(k/4)^{\frac{4}{k+4}}.
\end{equation*}
Take $t=\delta k$ and use $k^{6/(k+4)} \le e^6$ to get 
\begin{align*}
A_k &\lesssim  \frac{1}{\sqrt{\delta}} +  k  (\delta \CPIA(\rho))^{\frac{2}{k+4}} \Big(\delta\CPIA(\rho) + \frac{1}{2 k}\Big)^{\frac{k}{k+4}} (n-1)^{-\frac{2}{k+4}} \\
	&\lesssim  \frac{1}{\sqrt{\delta}} + k \big(\delta\CPIA(\rho) +1\big) (n-1)^{-\frac{2}{k+4}},
\end{align*}
Finally, optimize over $\delta$ to get

\begin{align*}
A_k &\lesssim  \Big(k \CPIA(\rho)(n-1)^{-\frac{2}{k+4}}\Big)^{1/3} + k (n-1)^{-\frac{2}{k+4}} \\
	&\lesssim  \big(1+\CPIA(\rho)^{1/3}\big)\big(k(n-1)^{-\frac{2}{k+4}}\big)^{1/3}.
\end{align*}
The last step is valid in the range $k(n-1)^{-2/(k+4)} \le 1$, and for the remaining range we can get away with the trivial bound $\BL \le 2$.
\end{proof}

\section{Sharpness proofs} \label{se:sharpness}

\subsection{Dependence on $\lambda$ in Theorem \ref{th:independent}} \label{se:kappa-dependence}

For a centered Gaussian $\rho$ with covariance matrix $\Sigma$, we can compute
\[
\HC(\rho) = \frac12 \log \frac{\prod_i \Sigma_{ii}}{\det \Sigma}.
\]
Consider $\rho$ being the centered Gaussian with covariance matrix $\Sigma=(1-\theta)I + \theta J$, for $0 < \theta < 1$, where $J$ is the all ones matrix. Note all diagonal entries are 1. The covariance matrix of the marginal $\rho_v$ for $|v|=k$ has the same form, and its determinant is $(1-\theta)^{k-1}(1+\theta (k-1))$. Thus, 
\[
\HC(\rho_v) = -\frac12 (k-1)\log(1-\theta) - \frac12 \log(1+\theta(k-1)) .
\]
The right-hand side is $\sim k(k-1)\theta^2/4$ for $\theta \to 0$ with $k$ fixed.
The largest and smallest eigenvalues of $\Sigma^{-1}$ are $b=1/(1-\theta)$ and $a=1/(1+\theta(n-1))$, respectively. A given $\lambda=b/a$ can be achieved by taking $\theta=(\lambda-1)/(n+\lambda-1)$, which yields
\[
\avg_{|v|=k} \HC(\rho_v) \sim \frac14 \frac{k(k-1)}{n^2}(\lambda-1)^2,
\]
for $k$ fixed as $n\to\infty$.
Thus, any coefficient $C(\lambda)$ satisfying
\[
\avg_{|v|=k} \HC(\rho_v) \le C(\lambda) \frac{k(k-1)}{n(n-1)}\HC(\rho),
\]
universally over $n$, $k$, and $\rho$, must satisfy
\begin{align*}
\frac14  (\lambda-1)^2 &\le C(\lambda) \lim_{n\to\infty} \Big(-\frac12 (n-1)\log(1-\theta) - \frac12 \log(1+\theta(n-1))\Big) \\
	&= \frac12 C(\lambda) \big(\lambda - 1  -  \log \lambda\big).
\end{align*}
In particular, $C(\lambda) > 1$ is necessary if $\lambda> 1$, and $C(\lambda)$ must grow at least like $\lambda/2$ as $\lambda \to\infty$.

\subsection{Proof of Proposition \ref{pr:sharpness-log}} \label{se:proof:pr:sharpness-log}

Fix $\delta > 2$ and let $k := \lceil \delta\log n / \log \log n\rceil$.
Let $q$ be the smallest power of 2 with $q \ge \sqrt{k}$, so that in particular $\sqrt{k} \le q < 2\sqrt{k}$. Set
\[
m = \lceil \log_q(n(q-1)+1)\rceil.
\]
Let $\FF_q$ be the finite field of order $q$. The number $m$ is chosen so that the number $(q^m-1)/(q-1)$ of one-dimensional subspaces of $\FF_q^m$ is at least $n$.
(There are $q^m-1$ non-zero vectors, and each one-dimensional subspace contains $q-1$ nonzero vectors.) We may thus choose $n$ nonzero vectors $a_1,\ldots,a_n\in\FF_q^m$ from distinct one-dimensional subspaces. Fix an arbitrary bijection $\iota : \FF_q \to \{0,1,\ldots,q-1\}$.

Let $U$ be drawn uniformly at random from the $q^m$ elements of $\FF_q^m$. Independently of $U$, let $V_1,\ldots,V_n$ be iid Unif$[0,1]$. Define
\[
X_i = \frac{1}{q}\big(\iota(a_i \cdot U) + V_i\big),
\]
and let $\rho$ be the law of $X=(X_1,\ldots,X_n)$.

We first claim that the $X_i$'s are pairwise independent Unif$[0,1]$. To see that they are uniform, note simply that $a_i \cdot U$ is uniform on $\FF_q$.
To see pairwise independence, note for $i \neq j$ that the vectors $a_i$ and $a_j$ are linearly independent. The map $\FF_q^m \ni u \mapsto (a_i\cdot u, a_j\cdot u) \in \FF_q^2$ is thus a surjection. Every point $(s,t) \in \FF_q^2$ has $q^{m-2}$ preimages under this map, and so
\[
\PP\big((a_i\cdot U, a_j\cdot U) = (s,t) \big) = \frac{q^{m-2}}{q^m} = q^{-2}.
\] 
This shows that $a_i \cdot U$ and $a_j \cdot U$ are independent, for $i \neq j$, and thus so are $X_i$ and $X_j$.

We next claim that $\CPIA(\rho)=1/12$. This will use the fact that each $X_i$ is $\mathrm{Unif}[0,1]$,  combined with the explicit calculation $\CPIA(\mathrm{Unif}[0,1])=1/12$. To see this, note that $\CPIA(\mathrm{Unif}[0,1])$  is the smallest $c$ such that $\Var(f(U)) \le c$ for all 1-Lipschitz $f$; here $U \sim \mathrm{Unif}[0,1]$. Note for any such $f$ that
\begin{align*}
\Var(f(U)) &\le \E[(f(U)-f(1/2))^2] \le \Var(U) = 1/12 .
\end{align*}
This becomes equality for $f(x)=x$, which establishes $\CPIA(\mathrm{Unif}[0,1])=1/12$. Now, turning back to $\rho$, note for any $f_1,\ldots,f_n$ that pairwise independence implies
\begin{align*}
\Var\Big(\sum_{i=1}^n f_i(X_i)\Big) = \sum_{i=1}^n \Var(f_i(X_i)) \le \frac{1}{12}\sum_{i=1}^n \|f'_i\|^2_{L^\infty(\R)}.
\end{align*}
This shows that $\CPIA(\rho)\le 1/12$, and equality again comes from taking linear functions. (Note that the statement of the proposition used $\CPIA(\rho)=1$ for simplicity, and this can be achieved by rescaling.)	

We next fix $v\subset [n]$ of size $k$ and choose a careful test function $f$ to bound the BL-distance from below. Define 
\begin{align*}
Q &:= \big\{(\iota(a_i\cdot u))_{i \in v} : u \in \FF_q^m\big\} \subset \{0,1,\ldots,q-1\}^k, \\
S &:= \bigcup_{z \in Q}\prod_{j=1}^k \Big[\frac{z_j}{q},\frac{z_j+1}{q}\Big].
\end{align*}
By construction, $\rho_v(S)=1$.
Define a 1-Lipschitz function $f:\R^k \to [0,1]$ by $f(x)=(1-d(x,S))_+$, where $d$ is the Euclidean distance.
This function satisfies $\int f\,d\rho_v = 1$. 
Note that each of the cubes in the definition of $S$ is contained in the ball of radius $\sqrt{k}/2q$ centered at $(2z+1)/2q$ where $1$ is the all-ones vector, and thus $f$ vanishes outside  of
\[
\bigcup_{z \in Q} B\Big(\frac{2z+1}{2q},1 + \frac{\sqrt{k}}{2q}\Big),
\]
where $B(y,r)$ denotes the ball of radius $r$ centered at $y$.
Because $\otimes_{i \in v}\rho_i$ is the uniform measure on $[0,1]^k$, and $|Q| \le q^m$, we deduce that 
\begin{align*}
\int f\,d\otimes_{i \in v}\rho_i = \int_{[0,1]^k} f(x)\,dx \le q^m\Big(\frac{2\pi e}{k}\Big)^{k/2}\Big(1 + \frac{\sqrt{k}}{2q}\Big)^k ,
\end{align*}
where $(2\pi e/k)^{k/2}$ is a standard upper bound on the volume of the unit ball in $\R^k$. Recalling that $\sqrt{k} \le q \le 2\sqrt{k}$, the definition of $m$ implies $q^m \le nq^2 \le 4nk$, and we bound the above by
\begin{align*}
\varepsilon_{n,k} :=  4n k^{1-\frac{k}{2}}(2\pi e)^{\frac{k}{2}} (3/2)^k.
\end{align*}
This construction worked for all $v$, and so
\[
\min_{|v|=k}\BL(\rho_v,\otimes_{i \in v}\rho_i) \ge 1- \varepsilon_{n,k}.
\]
Noting that
\[
\log\varepsilon_{n,k} \le \log n -\frac{k}{2}\log k + O(k),
\]
it follows quickly that $\varepsilon_{n,k} \to 0$ as $n\to\infty$ if $k \ge \delta\log n / \log \log n$ with $\delta>2$.

\begin{remark}
The ideas behind this construction originate in  the theory of \emph{orthogonal arrays}.
The $q^m \times n$ matrix $(a_i\cdot u)_{u \in \FF_q^m, i \in [n]}$ is called a (linear) orthogonal array, in design-theory terminology. Specifically, it is an array $\mathrm{OA}(q^m,n,q,2)$ of strength 2. This is sometimes known as the Rao-Hamming construction \cite[Section 3.4]{hedayat1999orthogonal}. 
Its use for us is as a way of producing as many pairwise-independent uniform $[q]$-valued random variables as possible from a sample space of size $q^m$.
In the language of coding theory, the set $\{(a_i \cdot u)_{i \in [n]} : u \in \FF_q^m\}$ is the $q$-ary simplex code, dual to the $q$-ary Hamming code. The precise construction of $X$, with the continuous uniform blurring, appeared in a paper of Owen \cite{owen1992orthogonal} in the context of Latin hypercube sampling.
\end{remark}

\bibliographystyle{amsplain}
\bibliography{8100bib}

\end{document}